\documentclass[11pt]{amsart}

\usepackage{amsfonts,amsmath}      
\usepackage{amssymb, latexsym}
\usepackage{amscd,amsthm}
\input xy
\xyoption{all}

\newcommand{\marginlabel}[1]%
  {\mbox{}\marginpar{\raggedleft\hspace{0pt}\bfseries\sf#1}}

\def\ZZ{{\mathbb Z}}

\def\QQ{{\mathbb Q}}
\def\CC{{\mathbb C}}

\def\OO{{\mathcal O}}
\def\R{\mathbf{R}}
\def\L{\mathbf{L}}
\def\SS{\mathcal{S}}
\def\LL{\mathcal{L}}
\def\D{\mathbf{D}}

\def\F{\mathcal{F}}
\def\E{\mathcal{E}}

\def\G{\mathcal{G}}
\def\H{\mathcal{H}}
\def\I{\mathcal{I}}
\def\J{\mathcal{J}}

\def\Pic0{{\rm Pic}^0(X)}

\theoremstyle{plain}

\newtheorem{theorem}{Theorem}[section]
\newtheorem{theoremalpha}{Theorem}
\newtheorem{corollaryalpha}[theoremalpha]{Corollary}

\newtheorem{proposition/example}[theorem]{Proposition/Example}
\newtheorem{proposition}[theorem]{Proposition}
\newtheorem{corollary}[theorem]{Corollary}
\newtheorem{lemma}[theorem]{Lemma}

\theoremstyle{definition}
\newtheorem{definition}[theorem]{Definition}
\newtheorem{remark}[theorem]{Remark}
\newtheorem{example}[theorem]{Example}

\newtheorem{conjecture/question}[theorem]{Conjecture/Question}

\newtheorem{remark/definition}[theorem]{Remark/Definition}
\newtheorem{definition/notation}[theorem]{Definition/Notation}

 \theoremstyle{remark}

\begin{document}

\title{GV-sheaves, Fourier-Mukai  transform, and Generic Vanishing}

\author[G. Pareschi]{Giuseppe Pareschi}
\address{Dipartamento di Matematica, Universit\`a di Roma, Tor Vergata, V.le della
Ricerca Scientifica, I-00133 Roma, Italy} \email{{\tt
pareschi@mat.uniroma2.it}}

\author[M. Popa]{Mihnea Popa}
\address{Department of Mathematics, University of Illinois at Chicago,
851 S. Morgan Street, Chicago, IL 60607, USA } \email{{\tt
mpopa@math.uic.edu}}

\thanks{MP was partially supported by the NSF grant DMS 0500985
and by an AMS Centennial Fellowship.}

\maketitle

\tableofcontents

\setlength{\parskip}{.1 in}

\markboth{G. PARESCHI and M. POPA} {\bf GV-sheaves, Fourier-Mukai
transform, and Generic Vanishing}

\section{Introduction}
In this paper we use homological techniques to establish a general
approach to generic vanishing theorems, a subject originated with the
pioneering work of Green-Lazarsfeld \cite{gl1} and \cite{gl2}.
Our work is inspired by a recent paper of Hacon
\cite{hacon}. Roughly speaking, we systematically investigate  --
in a general setting -- the relation between three concepts: \
(1) generic vanishing of (hyper-)cohomology groups of sheaves
(complexes) varying in a parameter space;  \ (2) vanishing of
cohomology \emph{sheaves} of Fourier-Mukai transforms; \  (3) a certain
vanishing condition for honest cohomology groups, related to
the vanishing of higher derived images in the spirit of the Grauert-Riemenschneider
theorem (which is a special case of this phenomenon). The relationship between these
concepts establishes a connection between Generic Vanishing theory
and the theory of Fourier-Mukai functors. One of the main points of
the paper is that, for projective varieties with mild singularities, these three
concepts are essentially equivalent and, more importantly, this equivalence
holds not just for derived equivalences, but in fact for any \emph{integral transform}.
This principle produces a number of new generic
vanishing results, which will be outlined below.

Let us briefly describe the three concepts above. For the sake of
brevity we do this only in the special case of smooth varieties and
locally free kernels (we refer to \S2,3 for a complete treatment, in
greater generality). Let $X$ and $Y$ be smooth projective
varieties, and let $P$ be a locally free sheaf on $X\times Y$. Here
$Y$ can be thought of as a moduli space of sheaves on $X$ and $P$ as the
universal sheaf. Generic vanishing theorems for cohomology groups
deal with the \emph{cohomological support loci}
$$V^i_P(\F)=\{y\in Y\>|\>h^i(\F\otimes P_y)>0\}$$
where $\F$ is a sheaf or, more generally, a complex on $X$ and
$P_y :=P_{|X\times \{y\}}$ is the locally free sheaf parametrized by $y$.
One can also consider the integral transform
$$\R \Phi_P : \D(X) \rightarrow \D (Y), \quad \R \Phi_P (\cdot) : = \R {p_Y}_* ( p_X^*(\cdot)
 \otimes P)$$ and the loci ${\rm Supp}(R^i \Phi_P \F)$ in $Y$. An
important point to note is that, although ${\rm Supp}(R^i \Phi_P
\F)$ is only contained in $V^i_P(\F)$, and in general not equal to it,
the sequences $\{\dim V^i_P(\F)\}_{i}$ and $\{\dim{\rm Supp}(R^i
\Phi_P \F)\}_i$ carry the same basic numerical information, in the sense that for any integer $k$
 \emph{the following conditions are equivalent}:

\noindent (a)
${\rm codim}_YV^i_P(\F)\geq i-k$ for all $i\ge 0$. \\
(b) ${\rm codim}_Y~{\rm Supp}(R^i \Phi_P \F)\geq i- k$  for all $ i
\ge 0.$

\noindent
(This follows from a standard argument using base change, see
Lemma \ref{equivalence} below.)  If these conditions are satisfied,
then $\F$ is said to \emph{satisfy Generic Vanishing with index $-k$}
with respect to $P$ or,  for easy reference, to be a
\emph{$GV_{-k}$-sheaf} (or \emph{$GV_{-k}$-object}).\footnote{In a previous 
version of this paper, we used the term \emph{$GV_k$-object} for what we 
now call a $GV_{-k}$-object. The reason for this change is that, as we will describe in upcoming work,
the notion of Generic Vanishing index is useful and can be studied for an arbitrary integer $k$, making 
more logical sense with this sign convention.}

Note that if $\F$ is a $GV_{-k}$-sheaf, then the cohomological
support loci $V^i_P(\F)$ are proper subvarieties for $i>k$. When
$k=0$ we omit the index, and simply speak of $GV$-sheaves (or
objects). An important example is given by the classical
Green-Lazarsfeld theorem \cite{gl1} which, in the above
terminology, can be stated as follows: \emph{let $X$ be a smooth
projective variety $X$, with Albanese map $a:X\rightarrow {\rm Alb} (X)$
and Poincar\' e line bundle $P$ on $X\times\Pic0$. Then $\omega_X$ is
$GV_{\dim a(X) -\dim X}$ with respect to $P$. In particular, if the
Albanese map of $X$ is generically finite, then $\omega_X$ is a
$GV$-sheaf.}

The second point -- sheaf vanishing -- concerns the vanishing of the higher
derived images $R^i\Phi_P \G$ for an object $\G$ in $\D(X)$. The question
becomes interesting when the locus $V^i_P(\G)$ is non-empty, since
otherwise this vanishing happens automatically.  Elementary
base change tells us that, for a potential connection with the $GV_{-k}$
condition  above,  the vanishing to look for has the following shape:
$$R^i\Phi_{P} \G = 0\>\>\> \hbox{for all \ } i < \dim X -k.$$
If this is the case,
we say that   \emph{$\G$ is $WIT_{\ge (\dim X-k)}$ with
respect to $P$}. The terminology is borrowed from Fourier-Mukai
theory, where an object $\G$ on $X$ is said \emph{to
satisfy the Weak Index Theorem (WIT) with index k} if
$R^i\Phi_P \G =0 $ for $i\ne k$. Continuing with the example of
the Poincar\'e bundle on $X\times\Pic0$, Hacon \cite{hacon} proved that $R^i\Phi_P \OO_X =0$
for $i<\dim a(X)$ i.e. -- in our terminology -- \emph{$\OO_X$ is
$WIT_{\ge \dim a(X)}$}. This was a conjecture of Green-Lazarsfeld
(cf. \cite{gl2} Problem 6.2;  see also \cite{pareschi} for a different argument).
In particular, if the Albanese map is generically finite, then $\R\Phi_P \OO_X$ is a
sheaf, concentrated in degree $\dim X$.

The third point expresses the vanishing of higher derived images  in terms of
the vanishing of a finite sequence of honest cohomology
\emph{groups}. A major step in this direction was made by Hacon for
the case when the $X$ is an abelian variety, $Y$ its dual, and $P$ is
a Poincar\'e line bundle. (In this case, by Mukai's theorem \cite{mukai}, the
functor $\R\Phi_P$ is an equivalence of categories.)
This was a key point in his proof in \cite{hacon} of the above mentioned
conjecture of Green-Lazarsfeld.

An essential point of the present work is that Hacon's argument can be suitably refined so that
it goes through for practically any integral transform $\R\Phi_P$, irrespective of whether
it is an equivalence, or even fully faithful.
To be precise, let us consider the analogous functor in the opposite direction
$$\R \Psi_P : \D(Y) \rightarrow \D (X),
~~\R \Psi_P (\cdot) : = \R {p_X}_* ( p_Y^*(\cdot) \otimes P).$$
For a sufficiently positive ample line bundle $A$ on $Y$, by Serre Vanishing we have
$\R\Phi_P(A^{-1})=R^g\Psi_P(A^{-1})[g]$. Moreover
$R^g\Psi_P(A^{-1})$ is locally free, and we denote it
$\widehat{A^{-1}}$. Following \cite{hacon}, one is lead to consider the cohomology groups
$H^i(X,\F\otimes \widehat{A^{-1}})$.  Denoting by $\R\Delta \F : = \R\mathcal Hom(\F,\omega_X)$
the (shifted) Grothendieck dual of $\F$, our general result is

\begin{theoremalpha}\label{f}
With the notation above, the following are equivalent:

\noindent
(1)  $\F$ is $GV_{-k}$.

\noindent
(2) $\R\Delta\F$ is $WIT_{\ge\dim X-k}$ with respect to
$\R\Phi_{P^\vee}:\D(X)\rightarrow \D(Y)$.

\noindent
(3) $H^i(X,\F\otimes\widehat{A^{-1}})=0$ for any $i>k$ and any sufficiently positive ample line
bundle $A$ on $Y$.
\end{theoremalpha}

We refer to Theorem  \ref{general_f} for the most general hypotheses
for Theorem \ref{f}: $X$ and $Y$ do not necessarily need to be
smooth, but rather just Cohen-Macaulay, and $P$
does not have to be a locally free sheaf, but rather any perfect object in
the bounded derived category of coherent sheaves. As already
observed in \cite{hacon}, condition (3) in Theorem \ref{f} is
extremely useful when $Y=\Pic0$ and $P$ is the Poincar\'e line
bundle on $X\times\Pic0$, since in this case $\widehat{A^{-1}}$ has
a very pleasant description: up to an \'etale cover of $X$, it it
the direct sum of copies of the pullback of an ample line bundle via
the Albanese map (cf. the proof of Theorem \ref{a} below). This
allows to reduce the Generic Vanishing conditions (1) or (2) to
classical vanishing theorems.

The implication $(2)\Rightarrow (1)$ is a rather standard application of
Grothendieck duality and of a basic fact on the support of Ext modules, 
well-known at least in the smooth case as a consequence of
the Auslander-Buchsbaum formula (cf. Proposition
\ref{WIT implies GV}). As mentioned above, the
equivalence $(2)\Leftrightarrow (3)$ was proved in
\cite{hacon} for $X$ and $Y$ dual abelian varieties and $P$ the
Poincar\'e line bundle. Hence the novel points of Theorem \ref{f}
are the implication $(1)\Rightarrow (2)$ and the equivalence
$(2)\Leftrightarrow (3)$ in the general setting of arbitrary
integral transforms. The latter is already important in
the well-studied case of the Poincar\'e line bundle $P$
on $X\times \Pic0$, for the following technical reason.
Deducing a result concerning the transform of $\F$ via
$\R\Phi_P:\D(X)\rightarrow \D(\Pic0)$ from a result involving the
derived equivalence $\R\Phi_{\mathcal{P}}:\D({\rm Alb}(X))\rightarrow \D(\Pic0)$,
using the Albanese map $a:X\rightarrow {\rm Alb} (X)$, requires
splitting and vanishing criteria for the object $\R a_* \F$ and its cohomologies.
Such criteria are available for the canonical bundle in the form of Koll\'ar's theorems on
higher direct images of dualizing sheaves \cite{kollar}, \cite{kollar2}. These however
require the use of Hodge theory, and are known not to hold in a more general setting,
for example for line bundles of the type $\omega_X\otimes L$ with $L$ nef.

We apply Theorem \ref{f} to deduce such a Kodaira-type
generalization of the Green-Lazarsfeld Generic Vanishing theorem to
line bundles of the form $\omega_X\otimes L$ with $L$ nef. When the
nef part is trivial, one recovers the above results of
Green-Lazarsfeld \cite{gl1} and Hacon \cite{hacon} -- in this case, although
the proof is just a variant of Hacon's, it has the following extra feature: Koll\'ar's
theorems on higher direct images of dualizing sheaves are not invoked, but rather just
Kodaira-Kawamata-Viehweg vanishing which, according to
Deligne-Illusie-Raynaud (\cite{di}, \cite{ev}), has an algebraic
proof via reduction to positive characteristic. Thus we provide a purely algebraic
proof of the Green-Lazarsfeld  Generic Vanishing theorem, answering a question of
Esnault-Viehweg (\cite{ev}, Remark 13.13(d)). To state the general theorem, we use the
following notation: for a $\QQ$-divisor $L$ on $X$, we define $\kappa_L$ to be $\kappa(L_{|F})$, the Iitaka dimension along the generic fiber of $a$, if $\kappa(L) \ge 0$, and $0$ if $\kappa(L) = -\infty$.

\begin{theoremalpha}\label{a}
Let $X$ be a smooth projective variety of dimension $d$ and Albanese
dimension $d-k$, and let $P$ be a Poincar\'e line bundle on $X\times
\Pic0$. Let $L$ be a line bundle and $D$ an effective $\QQ$-divisor
on $X$ such that $L-D$ is nef. Then $\omega_X\otimes L \otimes \J
(D)$ is a $GV_{-(k-\kappa_{L-D})}$-sheaf (with respect to $P$), where $\J (D)$ is the
multiplier ideal sheaf associated to $D$.  In particular,
if $L$ is a nef line bundle, then $\omega_X\otimes L$ is a $GV_{-(k-\kappa_L)}$-sheaf.
\end{theoremalpha}

\begin{corollaryalpha}\label{b}
Let $X$ be a smooth projective variety, and $L$ a nef line bundle on
$X$. Assume that either one of the following holds:

 \indent (1) $X$
is of maximal Albanese dimension.\hfill\break \indent (2) $\kappa(L)
\ge 0$ and $L_{|F}$ is big, where $F$ is the generic fiber of $a$.

\noindent Then  \ \ $H^i (\omega_X \otimes L \otimes P) = 0$ \ for
all $i
> 0$ and $P\in \Pic0$ general.   In particular
$\chi(\omega_X \otimes L) \ge 0$, and if strict inequality holds,
then $h^0 (\omega_X \otimes L) > 0$. The same is true if we replace
$L$ by $L\otimes \J(D)$, where $D$ is an effective $\QQ$-divisor on
$X$ and $L$ is a line bundle such that $L-D$ is nef.
\end{corollaryalpha}

\begin{corollaryalpha}\label{c}
If $X$ is a minimal smooth projective variety, then
$\omega_X^{\otimes m}$ is $GV_{-(k-\kappa_F)}$, where $\kappa_F$ is
equal to the Kodaira dimension $\kappa(F)$ of the generic fiber of
$a$ if $\kappa(F)\ge 0$, and $0$ otherwise.
\end{corollaryalpha}

Note that in the above results, a high Iitaka dimension along the
fibers compensates for the higher relative dimension of the Albanese
map. In the same vein, we obtain a number of new generic vanishing
results corresponding to standard vanishing theorems:

\noindent
$\bullet$ generic Koll\'ar-type vanishing criterion for higher direct image
sheaves of the form $R^if_*\omega_Y\otimes L$ with $L$ nef (Theorem \ref{direct}). \\
\noindent
$\bullet$ generic Nakano-type  vanishing criterion for bundles of
holomorphic forms (Theorem \ref{nakano_1}). \\
\noindent
$\bullet$ generic Le Potier, Griffiths and Sommese-type
vanishing for vector bundles (Theorem \ref{vector}).

In \S6  we give some first applications of these results to
birational geometry and linear series type questions, in the spirit
of \cite{kollar2} and \cite{el}. For instance, we show:

\noindent
$\bullet$ the multiplicativity of generic plurigenera under \'etale maps
of varieties whose Albanese map has generic fiber of
general type (Theorem \ref{plurigenera}).\\
\noindent
$\bullet$ the existence of sections, up to numerical equivalence,
for weak adjoint bundles on varieties of maximal Albanese dimension
(Theorem \ref{sections}).

\noindent
Here we exemplify only with a generalization of \cite{el} Theorem 3.

\begin{theoremalpha}\label{e}
Let $X$ be a smooth projective variety such that $V^0 (\omega_X)$ is
a non-empty proper subset of $\Pic0$. Then $a(X)$ is ruled by 
positive-dimensional subtori of
$A$.\footnote{Note that in case $X$ is of maximal Albanese dimension,
by generic vanishing the condition that $V^0(\omega_X)$ is a proper
subset of $\Pic0$ is equivalent to $\chi(\omega_X) = 0$.}
\end{theoremalpha}

We make some steps in the direction of higher rank Generic Vanishing
as well. We give an example of a generic vanishing result for
certain types of moduli spaces of sheaves on threefolds which are
Calabi-Yau fiber spaces (cf. Proposition \ref{cy} for the slightly
technical statement). The general such moduli spaces are pretty much
the only moduli spaces of sheaves on varieties of dimension higher
than two which seem to be quite well understood (due to work of
Bridgeland and Maciocia \cite{bm}, based also on work of Mukai).
A highly interesting point to be further understood here concerns
the description and the properties, in relevant examples, of the vector bundles
$\widehat{A^{-1}}$ used above. We also apply similar
techniques to moduli spaces of bundles on a smooth projective curve
to give a condition for a vector bundle to be a base point of any
generalized theta linear series, extending a criterion of Hein
\cite{hein} (cf. Corollary \ref{base_criterion_2}).

Finally, we mention that the results of this work have
found further applications of a different nature in the context of abelian
varieties in \cite{pp3} (which establishes a connection with $M$-regularity) and \cite{pp4}
(which uses Theorem \ref{f} to study subvarieties representing minimal cohomology classes
in a principally polarized abelian variety).\footnote{For the sake of
self-containedness,  an ad-hoc proof of Theorem
\ref{f} in the case of the Fourier-Mukai transform between dual
abelian varieties was given in \cite{pp4}, much simplified by the fact that
in this case one deals with an equivalence of categories.}

\noindent {\bf Acknowledgements.} The question whether there might
be a generic vanishing theorem for canonical plus nef line bundles
was posed to the first author independently by Ch. Hacon and M.
Reid. As noted above, we are clearly very much indebted to Hacon's
paper \cite{hacon}.  L. Ein has answered numerous questions and
provided interesting suggestions. We also thank J. Koll\'ar for
pointing out over-optimism in a general statement we had made about
surfaces, and to O. Debarre, D. Huybrechts and Ch. Schnell for useful
conversations and remarks. Finally, thanks are due to the 
referees for many valuable comments and corrections.

\section{Fourier-Mukai preliminaries}

We work over a field $k$ (in the applications essentially always over $\CC$). Given a 
variety $X$, we denote by $\D(X)$ bounded derived category of coherent sheaves on $X$.

Let $X$ and $Y$ be projective varieties over $k$, and $P$ a perfect object in 
$\D (X\times Y)$ (i.e. represented by a bounded complex of locally free sheaves of 
finite rank). This  gives as usual two Fourier-Mukai-type functors
$$\R \Phi_P : \D (X) \rightarrow \D(Y),
~~ \R \Phi_P (\cdot) : = \R {p_Y}_* ( p_X^*(\cdot) \underline\otimes
P),$$
$$\R \Psi_P : \D(Y) \rightarrow \D (X),
~~\R \Psi_P (\cdot) : = \R {p_X}_* ( p_Y^*(\cdot) \underline\otimes
P).$$

\noindent {\bf Projection formula and Leray isomorphism.} We will
use the same notation $H^{i} (\E)$ for the cohomology of a sheaf and
the hypercohomology of an object in the derived category. We recall
the following standard consequence of the projection formula and the
Leray isomorphism.

\begin{lemma}\label{projection}
For all objects $\E\in\D(X)$ and $\F\in\D(Y)$,
$$H^i(X,\E\underline\otimes \R\Psi_P \F)\cong H^i(Y,\R\Phi_P \E \underline\otimes \F).$$
\end{lemma}
\begin{proof}
By the projection formula (first and last isomorphism), and the
Leray equivalence (second and third isomorphism) we have:
\begin{eqnarray*}
{\bf R}\Gamma(X,\E\underline\otimes \R\Psi_P \F)&\cong
&\R\Gamma(X,\R{p_X}_* (p_X^*\E\underline\otimes
p_Y^*\F\underline\otimes P))\\ &\cong &{\bf R}\Gamma(X\times Y,
p_X^*\E\underline\otimes p_Y^*\F\underline\otimes P)\\ &\cong
&\R\Gamma(Y,\R{p_Y}_* (p_X^*\E\underline\otimes
p_Y^*\F\underline\otimes P))\\ &\cong &{\bf R}\Gamma(Y,\R\Phi_P \E
\underline\otimes \F).
\end{eqnarray*}
\end{proof}

Note that the statement works without any assumptions on the
singularities of $X$ and $Y$. We will use it instead of the more
common comparison of Ext groups for adjoint functors of
Fourier-Mukai functors on smooth varieties (cf. e.g. \cite{bo} Lemma
1.2).

\noindent {\bf Grothendieck duality.} We will frequently need to
apply Grothendieck duality to Fourier-Mukai functors. Given a
variety $Z$, for any object $\E$ in
$\D(Z)$ the derived dual of $\E$ is
$$\E^\vee : = \R \mathcal{H}om (\E, \OO_Z).$$
When $Z$ is Cohen-Macaulay and $n$-dimensional , we will also use the notation
$$\R \Delta \E: =  \R \mathcal{H}om (\E, \omega_Z),$$
so that the Grothendieck dualizing functor applied to $\E$ is $\R
\Delta \E [n]=  \R \mathcal{H}om (\E, \omega_Z [n])$.

Assume that $X$ and $Y$ are as above, with $X$ Cohen-Macaulay, 
and $P$ is an object in $\D(X\times Y)$. We use the notation introduced above.

\begin{lemma}\label{gd}
The Fourier-Mukai and duality functors satisfy the following
exchange formula:
$$(\R\Phi_P)^\vee \cong \R{\Phi}_{P^\vee} \circ \R \Delta_X [\dim X].$$
\end{lemma}
\begin{proof}
We have the following sequence of equivalences:
\begin{eqnarray*}
(\R\Phi_P (\cdot)) ^\vee &\cong & \R\mathcal{H}om( \R {p_Y}_* (
p_X^*(\cdot) \underline\otimes P), \OO_Y)\\
&\cong &\R {p_Y}_* ( \R\mathcal{H}om(p_X^*
(\cdot) \underline\otimes P,p_X^*\omega_X[\dim X]))\\
&\cong &\R {p_Y}_* ( \R\mathcal{H}om(p_X^*
(\cdot),p_X^*\omega_X[\dim X])\underline\otimes P^\vee)\\
&\cong &\R{p_Y}_*(p_X^* \R\Delta_X [\dim X] (\cdot)\underline\otimes
P^\vee )\\
&\cong &\R \Phi_{P^\vee}(\R\Delta_X(\cdot)[\dim X]).
\end{eqnarray*}
Besides basic operations allowed by the fact that we work with $X$
projective\footnote{What one needs is that the resolution
property for coherent sheaves be satisfied; cf. \cite{hartshorne}
Ch.II, \S5.}, the main point is Grothendieck Duality in the
second isomorphism. It works precisely as in the case of smooth
morphisms, given that in this case the relative dualizing sheaf for
$p_Y$ is $p_X^* \omega_X$ (cf. \cite{hartshorne} Ch.VII, \S4).
\end{proof}

\noindent {\bf Generalized WIT objects.} A key concept in
Fourier-Mukai theory is that of an object satisfying the Weak Index
Theorem, generalizing terminology introduced by Mukai \cite{mukai}
in the context of abelian varieties. We consider again $X$ and $Y$
projective varieties over $k$, and $P$ an object in $\D(X\times Y)$.

\begin{definition/notation}\label{wit_definition}
(1) An object $\F$ in $\D(X)$ is said to \emph{satisfy the Weak
Index Theorem with index $j$ ($WIT_j$ for short), with respect to
$P$}, if $R^i\Phi_P \F =0$ for $i\ne j$. In this case we denote
$\widehat \F=R^j\Phi_P \F $.

\noindent (2) More generally, we will say that $\F$ satisfies
$WIT_{\ge b}$ (or $WIT_{[a,b]}$ respectively) with respect to $P$,
if $R^i\Phi_P \F =0$ for $i<b$ (or for $i\not\in [a,b]$
respectively).
\end{definition/notation}

\begin{remark}[``Sufficiently positive"]\label{suff_positive}
In this paper we repeatedly use the notion of ``sufficiently positive" ample line bundle on $Y$ to mean, 
given any ample line bundle $L$, a power  $L^{\otimes m}$ with $m \gg0$. 
More precisely,  for a finite collection of coherent sheaves on $Y$, $m$ is sufficiently large so that  
by tensoring $L^{\otimes m}$ kills the higher cohomology of every sheaf in this collection.
\end{remark}

We have a basic cohomological criterion for
detecting $WIT$-type properties, as a  consequence of Lemma
\ref{projection} and Serre vanishing.

\begin{lemma}\label{cohomology}
Let $\F$ be an object in $\D(X)$. Then
$$R^i\Phi_P \F = 0 \iff H^i(X, \F\underline\otimes \R\Psi_P (A))=0$$
for any sufficiently positive ample line bundle $A$ on $Y$.
\end{lemma}
\begin{proof}
We use the spectral sequence
$$E^{ij}_2 := H^i(R^j\Phi_P \F \otimes A)\Rightarrow H^{i+j}(\R\Phi_P \F \otimes A).$$
By Serre vanishing, if $A$ is sufficiently ample, then
$H^i(R^j\Phi_P \F \otimes A)=0$ for $i>0$. Therefore $H^j(\R\Phi_P
\F \otimes A)\cong H^0(R^j\Phi_P \F \otimes A)$ for all $j$. Hence
$R^j\Phi_P\F$ vanishes if and only if $H^j(\R\Phi_P \F \otimes A)$
does. But, by Lemma \ref{projection}, $H^j(\R\Phi_P \F \otimes
A)\cong H^j(\F\underline\otimes \R\Psi_P A )$.
\end{proof}

\begin{corollary}\label{cohomology_cor}
Under the assumptions above, $\F$ satisfies $WIT_j$ with respect to
$P$ if and only if, for a sufficiently positive ample line bundle
$A$ on $Y$,
$$ H^i(X, \F\underline\otimes\R\Psi_P(A))=0 {\rm ~for~all~} i\ne j.$$
More generally, $\F$ satisfies $WIT_{\ge b}$ (or $WIT_{[a,b]}$
respectively) with respect to $P$ if and only if, for a sufficiently
positive ample line bundle $A$ on $Y$,
 $$ H^i(\F\underline\otimes\R\Psi_P(A))=0 {\rm ~for~all~} i<b ~ ({\rm or~} i\not\in [a,b]
 {\rm~respectively}).$$
\end{corollary}

\begin{example}[The graph of a morphism and Grauert-Riemenschneider]\label{graph}
The observation made in Lemma \ref{cohomology} is the generalization
of a well-known Leray spectral sequence method used in vanishing of
Grauert-Riemenschneider type. Let $f: X \rightarrow Y$ be a morphism
of projective varieties, and consider $P : = \OO_{\Gamma}$ as a
sheaf on $X\times Y$, where $\Gamma \subset X\times Y$ is the graph
of $f$. Hence $P$ induces the Fourier-Mukai functor $\R\Phi_P = \R
f_*$, and $\R\Psi_P$ is the adjoint $\L f^*$.

The criterion of Lemma \ref{cohomology} then says that for an object
$\F$ in $\D(X)$, we have
$$R^i f_* \F = 0 \iff H^i ( \F \otimes f^* A) = 0$$
for any $A$ sufficiently ample on $X$. This is of course well known,
and follows quickly from the Leray spectral sequence (cf.
\cite{positivity}, Lemma 4.3.10). For instance, if $X$ is smooth and $f$ has generic
fiber of dimension $k$, then $H^i (\omega_X \otimes L \otimes f^* A)
= 0$ for all $i >k$ and all $L$ nef on $X$, by an extension of
Kawamata-Viehweg vanishing (cf. Lemma \ref{iitaka} and the comments
before). This says that
$$R^i f_* (\omega_X\otimes L) = 0 , {\rm~for~all~} i > k,$$
which is a more general (but with identical proof) version of
Grauert-Riemenschneider vanishing.
\end{example}

\begin{example}[Abelian varieties]\label{abelian}
The statement of Lemma \ref{cohomology} is already of interest even
in the case of an abelian variety $X$ with respect to the classical
Fourier-Mukai functors $\R\hat{\mathcal{S}}: \D(X) \rightarrow
\D(\widehat X)$  and $\R \mathcal{S} : \D(\widehat X) \rightarrow
\D(X)$ (in Mukai's notation \cite{mukai}) given by a Poincar\'e line
bundle $\mathcal{P}$, i.e. in the present notation $\R
\Phi_{\mathcal{P}}$ and $\R \Psi_{\mathcal{P}}$. Since $A$ is
sufficiently positive, $\R \mathcal{S} A = \widehat A$ is a vector
bundle, and we simply have that $R^i \hat{\mathcal{S}} \F = 0$ iff
$H^i (\F \otimes \widehat A) = 0$. Note that $\widehat A$ is
negative, i.e. ${\widehat A}^\vee$ is ample (cf. \cite{mukai}
Proposition 3.11(1)).
\end{example}

\section{GV-objects}

In this section we introduce and study the notion of $GV$-object,
which is modeled on many geometrically interesting situations. We
will see that $WIT$ and $GV$-objects are intimately related,
essentially via duality. We use the notation of  the previous
section.

\begin{definition}
(1) An object $\F$ in $\D(X)$ is said to be a \emph{GV-object with
respect to $P$} if
$${\rm codim}~{\rm Supp}(R^i \Phi_P \F )\geq i {\rm ~for~all~} i\ge 0.$$

\noindent (2) More generally, for any integer $k \ge 0$, an object
$\F$ in $\D(X)$ is called a \emph{$GV_{-k}$-object with respect to $P$}
if
$${\rm codim}~{\rm Supp}(R^i \Phi_P \F)\geq i- k {\rm ~for~all~} i\ge 0.$$
Thus $GV = GV_0$. (This definition is short hand for saying that
$\F$ \emph{satisfies Generic Vanishing with index $k$ with respect
to $P$}.)

\noindent (3) If $\F$ is a sheaf, i.e. a complex concentrated in
degree zero,  and the conditions above hold, we call $\F$ a
\emph{$GV$-sheaf}, or more generally a \emph{$GV_{-k}$-sheaf}.
\end{definition}

In the next three subsections we establish the main technical
results of the paper, relating $WIT$ and  $GV$-objects. Together
with Lemma \ref{cohomology}, they will provide a cohomological
criterion for checking these properties.

\noindent {\bf GV versus  WIT.} All throughout we assume that $X$
and $Y$ are Cohen-Macaulay.  We set $d=\dim X$ and $g=\dim Y$. 

\begin{proposition}[$GV$ implies $WIT$]\label{GV implies WIT}
Let $\F$ be an object in $\D(X)$ which is $GV_{-k}$  with respect to
$P$. Then $\R\Delta \F$ is $WIT_{\ge (d-k)}$ with respect to $P^\vee \underline\otimes p_Y^* \omega_Y$.
\end{proposition}
\begin{proof}
Denote $Q := P^\vee \underline\otimes p_Y^* \omega_Y$.
Let $A$ be a sufficiently positive ample line bundle on $Y$. By
Lemma \ref{cohomology}, it is enough to prove that
$$H^i(X,\R\Delta \F \underline\otimes \R\Psi_Q A)=0 {\rm ~for~all~} i < d-k.$$
By Grothendieck-Serre duality, this is equivalent to
$$H^j(X, \F\underline\otimes (\R\Psi_Q A)^\vee)=0 {\rm~for~all~} j > k.$$
By Lemma \ref{gd} we have $(\R\Psi_Q A)^\vee\cong \R\Psi_{Q^\vee}(\R
\Delta A)[g] \cong \R\Psi_P (A^{-1})[g]$, where the second
isomorphism follows by a simple calculation. In other words, what we
need to show is
$$H^{j +g} (X, \F \underline\otimes \R\Psi_P (A^{-1}))=0 {\rm ~for~all~} j > k.$$
This in turn is equivalent by Lemma \ref{projection} to
$$H^l (Y, \R \Phi_P \F \otimes A^{-1}) = 0 {\rm ~for~all~} l >  g + k.$$
Now on $Y$ we have a spectral sequence
$$E^{pq}_2 := H^p  (R^q \Phi_P \F \otimes A^{-1})\Rightarrow
H^{p+q} (\R \Phi_P \F \otimes A^{-1}).$$ Since $\F$ is $GV_{-k}$ with
respect to $P$, whenever $l = p+q > g+k$ we have that ${\rm dim}
~{\rm Supp} (R^q \Phi_P \F) < p$, and therefore  $E^{pq}_2 = 0$.
This implies precisely what we want.
\end{proof}

\begin{proposition}[$WIT$ implies $GV$]\label{WIT implies GV}
Let $\F$ be an object in $\D(X)$ which satisfies $WIT_{\ge (d-k)}$
with respect to $P \underline\otimes p_Y^* \omega_Y$. 
Then $\R\Delta \F$ is $GV_{-k}$ with respect to $P^\vee$.
\end{proposition}
\begin{proof}
Grothendieck duality (Lemma \ref{gd}) gives
$$\R\Phi_{P^\vee}\bigl(\R\Delta \F \bigr) \cong (\R\Phi_P \F [d])^\vee
\cong \R \Delta (\R\Phi_{P\underline \otimes p_Y^* \omega_Y} \F [d]).$$
By assumption, $\R\Phi_{P\underline \otimes p_Y^* \omega_Y} \F [d]$ is an object whose cohomologies $R^j$
are supported in degrees at least $-k$. 
The spectral sequence
associated to the composition of two functors implies in this case
that there exists a spectral sequence
$$E^{ij}_2 := \E xt^{i+j} (R^j, \omega_Y)\Rightarrow R^i\Phi_{P^\vee} (\R \Delta \F).$$
One can now use general facts on the support of Ext-sheaves. More
precisely, since $Y$ is Cohen-Macaulay, we know that
$${\rm codim} ~{\rm Supp}(\mathcal{E}xt^{i + j}(R^j, \omega_Y)) \ge i + j$$
for all $i$ and $j$. (This is better known as an application of the Auslander-Buchsbaum  
formula when the $R^j$ have finite projective dimension, but holds in general by 
e.g. \cite{bh} Corollary 3.5.11(c).) Since the only non-zero
$R^j$-sheaves are for $j \ge -k$, we have that the codimension of
the support of  every $E_{\infty}$ term is at least $i -k$. This
implies immediately that ${\rm codim}~{\rm Supp}(R^i \Phi_{P^\vee} (\R \Delta \F))
\geq i-k$, for all $i\ge 0$.
\end{proof}

\begin{corollary}\label{hacon implies GV}
Let $\F$ be an object in $\D(X)$ such that
$$H^i(\F \underline\otimes \R\Psi_{P [g]} (A^{-1}))=0 {\rm ~for~all~} i > k,$$
for any sufficiently positive ample line bundle $A$ on $Y$. Then
$\F$ is $GV_{-k}$ with respect to $P$.
\end{corollary}
\begin{proof}
By Proposition \ref{WIT implies GV}, it is enough to prove that
$\R\Delta \F$ is $WIT_{\ge (d-k)}$ with respect to $Q : = P^\vee \underline\otimes p_Y^* \omega_Y$. 
Using
Corollary \ref{cohomology_cor} and Grothendieck-Serre duality, this
follows as soon as
$$H^i(X,\F\otimes ({\bf R}\Psi_{Q} A)^\vee)=0 {\rm ~for~ all~} i > k.$$
We are left with noting that, by Lemma \ref{gd},
$$({\bf R} \Psi_{Q} A )^\vee\cong
\R\Psi_{P \underline\otimes p_Y^*\omega_Y^\vee  [g]}(A^{-1} \otimes \omega_Y) 
\cong \R\Psi_{P [g]} (A^{-1}),$$
the last isomorphism being due to the Projection Formula and the fact that $\omega_Y \underline
\otimes \omega_Y^\vee \cong \OO_Y$.
\end{proof}

\noindent {\bf Cohomological support loci.} 
Generic Vanishing
conditions were originally given in terms of cohomological support
loci, so it is natural to compare the definition of $GV_{-k}$-sheaves
with the condition that the $i$-th cohomological support locus of
$\F$ has codimension $\ge i-k$. For any $y\in Y$ we denote $P_y= \L
i_{y}^* P$ the object in $\D(X)$, where $i_{y}: X\times \{y\}
\hookrightarrow X \times Y$ is the inclusion.

\begin{definition}
Given an object $\F$ in $\D(X)$, the $i$-th cohomological support
locus of $\F$ with respect to $P$ is
$$V^i_P(\F) : = \{y\in Y\>|\> \dim H^i(X,\F \underline\otimes P_y)>0\}.$$
\end{definition}

\begin{lemma}\label{equivalence}
The following conditions are equivalent:\\
(1) $\F$ is a $GV_{-k}$-object with respect to $P$. \\
(2) ${\rm codim}_YV^i_P(\F)\geq i-k$ for all $i$.
\end{lemma}
\begin{proof}\footnote{Cf. also \cite{hacon} Corollary 3.2 and \cite{pareschi} Corollary 2.}
The statement follows from the theorem on cohomology and base change
and its consequences.\footnote{We recall that, although most
commonly stated for cohomology of coherent sheaves (see e.g. the
main Theorem of \cite{mumford}, \S 5 p.46), base change  -- hence
also its corollaries, as \cite{mumford} Corollary 3, \S 5 -- works
more generally for hypercohomology of bounded complexes 
(\cite{ega3} 7.7, especially 7.7.4, and Remarque 7.7.12(ii)). }
Since by cohomology and base change ${\rm Supp}(R^i \Phi_P
\F)\subseteq V^i_P(\F)$, it is enough to prove that (1) implies (2).
The proof is by descending induction on $i$. There certainly exists
an integer  $s$, so that $H^{j}(\F\underline\otimes P_y)=0$ for any
$j>s$ and for any $y\in Y$. Then, by base change, ${\rm Supp}(R^s
\Phi_P \F)=V^s_P(\F)$. The induction step is as follows: assume that
there is a component $\bar V$ of $V^i_P(\F)$ of codimension less
than $i-k$. Since (1) holds, the generic point of $\bar V$ cannot be
contained in ${\rm Supp}(R^i \Phi_P \F)$ and so, again by base
change, we have that $\bar V\subset V^{i+1}_P(\F)$. This implies
that ${\rm codim}_Y V^{i+1}_P(\F)< i - k$, which contradicts the
inductive hypothesis.
\end{proof}

Here we will only use this in the standard setting where $P$ and
$\F$ are sheaves, with $P$ locally free. In this case $P_y$ is
just the restriction of $P$ to $X\times\{y\}$ and $V^i_P(\F)$ is
simply the locus where the sheaf cohomology $H^i (\F \otimes P_y)$
is non-zero.

\noindent {\bf The main technical result.} 
One can put together the sequence of results above in order
to obtain the main technical result, implying Theorem \ref{f}
in the Introduction.

\begin{theorem}\label{general_f}
Let $X$ and $Y$ be projective Cohen-Macaulay varieties, of  dimensions 
$d$ and $g$ respectively, and let $P$ be a perfect object in $\D(X\times Y)$. 
Let $\F$ be an object in $\D(X)$. The following conditions are equivalent:

\noindent (1) $\F$ is a $GV_{-k}$-object with respect to $P$.

\noindent (2) $H^i(\F\underline\otimes \R\Psi_{P [g]} (A^{-1}))=0$
for $i > k$ and any sufficiently positive ample line bundle $A$ on $Y$.

\noindent (3) $R^i\Phi_{P^\vee \underline \otimes p_Y^* \omega_Y} (\R \Delta \F) = 0$ for all $i <
d-k$ (i.e. $\R \Delta \F$ satisfies $WIT_{\ge (d-k)}$ with respect
to $P^\vee \underline \otimes p_Y^* \omega_Y$).

\noindent (4) ${\rm codim}_YV^i_P(\F)\geq i-k$ for all $i$.
\end{theorem}

\begin{proof}
Everything follows by putting together Proposition \ref{GV implies
WIT}, Proposition \ref{WIT implies GV}, Corollary \ref{hacon implies
GV} and Lemma \ref{equivalence}. 
\end{proof}

\begin{remark}
To make the transition to the statement of Theorem \ref{f} in the Introduction, simply note
that we have 
$$\R\Phi_{P^\vee \underline \otimes p_Y^* \omega_Y} (\R \Delta \F) \cong 
\R\Phi_{P^\vee} (\R \Delta \F) \underline \otimes \omega_Y.$$ 
In case $Y$ is Gorenstein, $\omega_Y$ is a line bundle so (3) above is equivalent 
to $\R \Delta \F$ being $WIT_{\ge (d-k)}$ with respect to $P^\vee$.
\end{remark}

\begin{example}[The graph of a morphism II]\label{graph2}
As in Example \ref{graph}, let $f: X \rightarrow Y$ be a morphism of
projective varieties, and consider $P : = \OO_{\Gamma}$ as a sheaf
on $X\times Y$, where $\Gamma \subset X\times Y$ is the graph of
$f$. We have $\R\Phi_P = \R f_*$ and $\R\Psi_P = \L f^*$. Assuming
that $X$ and $Y$ are Cohen-Macaulay, the interpretation
of Theorem \ref{general_f} in this case is that an object $\F$ in
$\D(X)$ is $GV_{-k}$ with respect to $P$ if and only if $H^{i}(\F
\otimes f^* (A^{-1})) = 0$ for all $i >  g+ k$ and any $A$
sufficiently positive on $Y$. In other words, in analogy with
Example \ref{graph} we have the following, presumably folklore,
consequence:

\begin{corollary}\label{push_forward}
If $f$ has generic fiber of dimension $k$, then for any object $\F$
in $\D(X)$:
$${\rm codim}~{\rm Supp}(R^i f_* \F) \ge i-k \iff H^i ( \F \otimes f^* (A^{-1})) = 0, ~\forall i > g + k.$$
For example ${\rm codim}~{\rm Supp}(R^i f_* \OO_X) \ge i-k, {\rm
~for~all~} i$.
\end{corollary}
\end{example}

\medskip

In most instances in which Theorem \ref{general_f} is applied, due
to geometrically restrictive assumptions on $\F$ and $P$, it is also
the case that $R^i \Phi_{P^\vee} (\R\Delta \F) = 0$ for $i >d$.

\begin{corollary}\label{special}
If in Theorem \ref{general_f} we assume in addition that the kernel
$P$ and $\F$ are  sheaves, with $P$ locally free, then $\F$ being $GV_{-k}$ with
respect to $P$ is equivalent to

\noindent (2$^\prime$) $H^i(\F \otimes \R\Psi_{P [g]}
(A^{-1}))=0$ for $i \not\in[0, k]$, and for any sufficiently
positive ample line bundle $A$ on $Y$.

\noindent (3$^\prime$) $R^i\Phi_{P^\vee \otimes  p_Y^* \omega_Y} (\R \Delta \F) = 0$ for all
$i \not\in [d-k,d]$ (i.e. $\R \Delta \F$ satisfies $WIT_{[d-k, d]}$
with respect to $P^\vee \otimes p_Y^* \omega_Y$).
\end{corollary}
\begin{proof}
By Lemma \ref{cohomology}, the only thing we need to note is that
under the extra assumptions we have $R^i \Phi_{P^\vee \otimes p_Y^* \omega_Y} (\R\Delta \F)
= 0$ for $i>d$. Indeed, since $\F$ is a sheaf and $P$ is locally free, $V^{d-i}_P(\F)=V^i_{P^\vee}(\R\Delta\F)$ (Serre duality) is empty for $i>d$. But $V^i_{P^\vee}(\R\Delta\F)=V^i_{P^\vee\otimes p_Y^*\omega_Y}(\R\Delta\F)$,  so the assertion follows by base change.
\end{proof}

For convenience, it is worth stating the result of Theorem
\ref{general_f} in the most important special case, namely the
relationship between $GV= GV_0$ and $WIT_d$ under these geometric
assumptions.

\begin{corollary}\label{k=0}
Under the assumptions of Theorem \ref{general_f}, if 
$P$ and $\F$ are sheaves, and $P$ is locally free,
the following conditions are  equivalent:

\noindent (1) $\F$ is $GV$ with respect to $P$.

\noindent (2) $H^i(\F \otimes \R\Psi_{P [g]} (A^{-1}))=0$
for $i \neq 0$, for a sufficiently positive ample line bundle $A$ on
$Y$.

\noindent (3) $\R\Delta \F$ satisfies $WIT_d$ with respect to
$P^\vee \otimes p_Y^* \omega_Y$.

\noindent (4) ${\rm codim}_YV^i_P(\F)\geq i$ for all $i$.
\end{corollary}

\begin{remark}\label{gv-duality}
We emphasize that, as it follows from the proof (of Proposition
\ref{WIT implies GV}), if the equivalent conditions of Corollary
\ref{k=0} hold, then
$$R^i\Phi_P \F \cong \mathcal{E}xt^i(\widehat{\R\Delta \F},\omega_Y),$$
where the Fourier-Mukai hat is taken with respect to $P^\vee$ (cf.
Definition/Notation \ref{wit_definition}).
\end{remark}

Another remark which is very useful in applications (see e.g.
\cite{el}) is the following (cf. also \cite{hacon} Corollary 3.2(1)
and \cite{pareschi} Corollary 2):

\begin{proposition}\label{general_inclusions}
Assume that $P$ is locally free. Let $\F$ be a $GV_{-k}$-sheaf with
respect to $P$. Then
$$V^d_P(\F) \subseteq \ldots \subseteq V^{k+1}_P (\F)\subseteq V^{k}_P (\F).$$
\end{proposition}
\begin{proof}
By Grothendieck-Serre duality we have for all $i$ and all $y \in Y$
that
$$H^i (\F \otimes P_y) \cong H^{d-i} (\R\Delta \F \otimes P_y^\vee)^\vee.$$
If this is $0$, then by base change (cf.  e.g. \cite{mumford} \S5,
Corollary 2, and also the comments in the proof of Lemma \ref{equivalence}) the natural
homomorphism
$$R^{d-i -1} \Phi_{P^\vee} (\R\Delta \F)\otimes k(y) \longrightarrow
H^{d-i-1} (\R\Delta \F \otimes P_y^\vee)$$ is an isomorphism. Since
$\F$ is $GV_{-k}$ with respect to $P$, Theorem \ref{general_f} and 
the Projection Formula imply
that $R^{d-i -1} \Phi_{P^\vee} (\R\Delta \F)$ is $0$ for $i \ge k$,
so again by duality we get that $H^{i+1} (\F \otimes P_y)= 0$.
\end{proof}

\noindent {\bf The case when $V^0$ is a proper subvariety.} Assume
that $\R\Phi_P$ is a Fourier-Mukai functor between $X$
and $Y$ projective Cohen-Macaulay, with $P$ a locally free sheaf.
When $V^0(\F)$ is a proper subvariety of $Y$, one has strong
consequences, useful in geometric applications (cf. \S7). This
generalizes results of Ein-Lazarsfeld for Albanese maps (cf.
\cite{el} \S1 and \S2).

\begin{proposition}\label{component}
Assume that $\F$ is a  sheaf on $X$ which is $GV$ with
respect to $P$, and let $W$ be a component of $V^0 (\F)$ of
codimension $p$. If $W$ is either isolated or of maximal dimension
among the components of $V^0(\F)$, then ${\rm dim}~X \ge p$, and $W$
is also a component of $V^p (\F)$.
\end{proposition}
\begin{proof}
By Grothendieck-Serre duality, the sheaf $\G : =
\widehat{\R\Delta\F}$ on $Y$ (as in Remark \ref{gv-duality})
has support $V^0 (\F)$. Denote $\tau:=
\G_{|W}$ and consider the exact sequence given by restriction to
$W$:
$$0 \rightarrow \H \rightarrow \G \rightarrow \tau \rightarrow 0.$$
Assuming that $W$ is of maximal dimension among the components of $V^0 (\F)$, 
the codimension of the support of $\H$ is at least $p$, which
implies that $\mathcal{E}xt^i (\H, \omega_Y) = 0$, for all $i < p$. As
a consequence, for $i\le p$ we have an inclusion $\mathcal{E}xt^i
(\tau, \omega_Y) \hookrightarrow \mathcal{E}xt^i (\G, \omega_Y)$. This last assertion 
is obviously also true if $W$ is an isolated component of $V^0 (\F)$.  On the
other hand, it is standard that the support of
$\mathcal{E}xt^p(\tau,\omega_Y)$ is $W$.

By Remark \ref{gv-duality} we know that $\E xt^{i} (\G, \omega_Y) \cong
R^i \Phi_P \F$, which is $0$ for $i >d$. Combined with the above,
this gives $p\le \dim X$. In addition $V^p(\F)$ contains the support
of $\mathcal{E}xt^p(\G,\omega_Y)$, which must contain the support of
$\mathcal{E}xt^p(\tau,\omega_Y)$, i.e $W$. But since $\F$ is $GV$, $V^p
(\F)$ has codimension at least $p$, so $W$ must then be one of its
components.
\end{proof}

\begin{corollary}\label{isolated}
Assume that $\F$ is a sheaf on $X$ which is $GV$ with
respect to $P$, such that $V^0 (\F)$ has an isolated point. Then
${\rm dim}~X \ge {\rm dim}~Y$.
\end{corollary}

\section{Examples of $GV$-objects}

We have seen some basic examples related to morphisms in \ref{graph}
and \ref{graph2}. In what follows we present a few other, more
interesting, examples of $GV$-objects. With the exception of (5),
all are in the following context: $X$ is a smooth projective variety
of dimension $d$, with Albanese map $a: X \rightarrow A$, $P$ is a
Poincar\'e  bundle on $X \times \widehat A$, and $\R\Phi_P : \D(X)
\rightarrow \D(\widehat A)$ is the induced Fourier-Mukai functor.
The example in (5) is related to \S8.

\noindent (1) {\bf  The Green-Lazarsfeld Generic Vanishing theorem.}
The main theorem of \cite{gl1} says in the present language that if
$a$ has generic fiber of dimension $k$, then $\omega_X$ is a
$GV_{-k}$-sheaf with respect to $P$.  Green and Lazarsfeld also
conjectured that if $a$ is generically finite ($k = 0$), then $R^i
\Phi_P \OO_X = 0$ unless $i=d$, i.e. that $\OO_X$ satisfies $WIT_d$
with respect to $P$. This, and actually that in the general case
$\OO_X$ satisfies $WIT_{[d -k,d]}$, was proved by Hacon \cite{hacon}
and Pareschi \cite{pareschi}. Theorem \ref{general_f} and Theorem
\ref{a} imply that this is in fact equivalent to the theorem in
\cite{gl1}.

\noindent (2) {\bf Line bundles on curves.} Let $X = C$, a smooth
projective curve of genus $g$, so that $A \cong \widehat A\cong
J(C)$, and $a$ is an Abel-Jacobi map. If $L$ is a line bundle on
$C$, then $L$ is $GV$ with respect to $P$ if and only if ${\rm
deg}(L) \ge g-1$. Dually, $L$ satisfies $WIT_1$ if and only if ${\rm
deg}(L) \le g-1$. Every line bundle on $C$ is $GV_{-1}$.

\noindent (3) {\bf Line bundles on symmetric products and a result
of Polishchuk.} A natural generalization of the previous example is
as follows. For $1 \le d \le g-1$, let $X = C_d$ be the $d$-th
symmetric product of $C$. Its Albanese variety is $J(C)$ and the
Albanese map is an abelian sum mapping $a:C_d\rightarrow J(C)$. Let
also $\pi_d:C^d\rightarrow C_d$ be the natural projection. To a line
bundle $L$ on $C$, one can attach canonically a line bundle $F_d(L)$
on $C_d$ such that $\pi_d^*F_d(L) \cong L^{\boxtimes d}$.

We claim that $\omega_{C_d}\otimes (F_d(L))^{-1}$ is $GV$ with
respect to $P^\vee$ (so also $P$) if and only if $\deg L\le g-d$. In
order to prove this, we recall some standard facts (cf. e.g.
\cite{izadi} Appendix 3.1). In the first place, if $\xi \in {\rm
Pic}^0 (C)$ and $P_{\xi}$ is the corresponding line bundle on
$J(C)$,  then $F_d(L)\otimes a^*P_\xi\cong (L\otimes
P_\xi)^{\boxtimes d}$. Moreover
$\pi^*\omega_{C_d}\cong\omega_C^{\boxtimes d}(-\Delta)$, where
$\Delta$ is the sum of the big diagonals in $C^d$. Therefore
$$\pi_d^*(\omega_{C_d}\otimes (F_d(L)\otimes
a^*P_\xi)^{-1})\cong (\omega_C\otimes (L\otimes
P_\xi)^{-1})^{\boxtimes d}(-\Delta).$$ The cohomology of
$\omega_{C_d}\otimes (F_d(L)\otimes a^*P_\xi)^{-1}$ is the
skew-symmetric part of the cohomology of $(\omega_C\otimes (L\otimes
\xi)^{-1})^{\boxtimes d}$ with the respect to the action of the
symmetric group, so
$$H^i(C_d,\omega_{C_d}\otimes (F_d(L)\otimes
a^*P_\xi)^{-1})\cong S^i(H^1(C, \omega_C\otimes
(L\otimes\xi)^{-1}))\otimes \Lambda^{d-i} (H^0(C, \omega_C\otimes
(L\otimes\xi)^{-1})).$$ A simple calculation shows that, for any
$i$:
$${\rm codim} V^i_{P^\vee}(\omega_{C_d}\otimes
(F_d(L))^{-1})\ge i \iff \deg L\le g-d.$$

Although we mainly focus here on the $GV$ notion, Theorem
\ref{general_f} can also be used dually to check that an object
satisfies $WIT$, which is sometimes harder to prove. For example,
using (1) $\Rightarrow$ (3), the calculation above gives a quick
proof of a result of Polishchuk which essentially says:

\begin{corollary}[cf. \cite{polishchuk2} Theorem 0.2]
The line bundle $F_d(L)$ satisfies $WIT_d$ with respect to $P$ if
and only if $\deg L\le g-d$.
\end{corollary}

\noindent (4) {\bf Ideal sheaves of subvarieties in ppav's. }
Consider now $X = J(C)$, the Jacobian of a smooth projective curve
of genus $g$, with principal polarization $\Theta$.  For each $d =
1, \ldots, g-1$, denote by $W_d$ the variety of special divisors of
degree $d$ in $J(C)$. In \cite{pp1} Proposition 4.4, we proved that
the sheaves $\OO_{W_d} (\Theta)$ are $M$-regular and
$h^0(\OO_{W_d}(\Theta) \otimes \alpha) = 1$ for $\alpha \in \Pic0$
general. Using the twists of the standard exact sequence
$$0\longrightarrow \I_{W_d}(\Theta)\longrightarrow \OO_{J(C)}(\Theta)\longrightarrow \OO_{W_d}(\Theta) \longrightarrow 0$$
by $\alpha\in \Pic0$, and Lemma \ref{equivalence}, we deduce easily
that $\I_{W_d}(\Theta)$ are $GV$-sheaves on $J(C)$. In the paper
\cite{pp4} we use Theorem \ref{general_f} to go in the opposite
direction: we show, among other things, that if $\I_Y(\Theta)$ is a
$GV$-sheaf then $Y$ represents a \emph{minimal} class. This, in
combination with a conjecture of Beauville-Debarre-Ran, means that the existence
of subvarieties $Y$ in a ppav $(A, \Theta)$ such that $\I_Y
(\Theta)$ is $GV$ should always characterize Jacobians, with the
exception of intermediate Jacobians of cubic threefolds.

\noindent (5) {\bf Stable sheaves on surfaces.} Consider $X$ to be a
complex abelian or $K3$ surface.  For a coherent sheaf $E$ on $X$,
the  Mukai vector of $E$ is
$$v(E) : = {\rm rk} (E) + c_1(E) + (\chi (E) -  \epsilon \cdot rk(E) ) [X] \in H^{ev}(X, \ZZ),$$
where $\epsilon$ is $0$ if $X$ is abelian and $1$ if $X$ is $K3$.
Given a polarization $H$ on $X$ and a vector $v \in H^{ev}(X, \ZZ)$,
we can consider the moduli space $M_H (v)$ of sheaves $E$ with $v(E)
= v$, stable with respect to $H$. If the Mukai vector $v$ is
primitive and isotropic, and $H$ is general, the moduli space is $M
= M_H (v)$ is smooth, projective and fine, and it is in fact again
an abelian or $K3$ surface (cf. e.g. \cite{yoshioka1}). The
universal object $\E$ on $X \times M$ induces an equivalence of
derived categories $\R \Phi_{\E} : \D(X) \rightarrow \D(M)$.
Yoshioka gives many examples of situations in which in our language
the $GV$ (or, dually, the $WIT$)  property is satisfied by $\OO_X$
($\cong \omega_X$) with respect to $\E$. Here we give just a very
brief sampling. Some of these results will be used in Proposition
\ref{cy}.

Consider for example $(X, H)$ to be a polarized $K3$ surface such
that ${\rm Pic}(X) = \ZZ\cdot H$, with $H^2 = 2n$. Let $k > 0$ be an
integer, and consider $v = k^2n + kH + [X]$. This is a primitive
isotropic Mukai vector. Assume in addition that $kH$ is very ample.
It is shown in \cite{yoshioka3} Lemma 2.4 that under these
assumptions $WIT_2$ holds for $\OO_X$ with respect to the kernel
$\E^\vee$.\footnote{It is also shown in \emph{loc. cit.} \S2 that in
fact $X \cong M$.}  By Corollary \ref{k=0}, this is equivalent to
the fact that $\OO_X$ is a $GV$-sheaf with respect to $\E$.

There are similar examples when $(X, H)$ is a polarized abelian
surface with ${\rm Pic}(X) = \ZZ\cdot H$. Write $H^2 = 2r_0 k$, with
$(r_0,k) = 1$. Consider the Mukai vector $v_0 := r_0 + c_1(H) +
k[X]$, which is primitive and isotropic, so $M = M_H(v_0)$ is again
an abelian surface, and there exists a universal object $\E$ on $X
\times M$. It is proved in \cite{yoshioka2}, Theorem 2.3 and the
preceding remark, that $\OO_X$ (among many other examples) satisfies
$WIT_2$ with respect to $\E^\vee$, which as above means that $\OO_X$
is $GV$ with respect to $\E$.

If $(X,H)$ is a polarized abelian surface, one can consider the
behavior of individual stable bundles with respect to the usual
Fourier-Mukai functor $\R\hat{\mathcal{S}} : \D(X) \rightarrow
\D(\widehat X)$ as well. Assume for example that $NS(X) = \ZZ\cdot
H$, and consider the Mukai vector $v = r + c_1(L) + a [X]$, where $a
> 0$. It is shown in \cite{yoshioka2} Proposition 1.1 that if $E$ is
a stable bundle with respect to $H$, with Mukai vector $v$, then $E$
satisfies $WIT_2$ with respect to the dual Poincar\'e bundle
$\mathcal{P}^\vee$, or in other words that $E^\vee$ is a $GV$-sheaf
with respect to $\R\hat{\mathcal{S}}$.

\section{Generic vanishing theorems}

In this section we give the main applications of the material in
\S3, namely Generic Vanishing theorems related to the Picard
variety. The main statement is phrased in the context of multiplier ideal
sheaves. It generalizes results in \cite{gl1}, \cite{hacon}, and
\cite{pareschi}, and it contains as a special case generic vanishing
for adjoint bundles of the form $K_X + L$ with $L$ nef. We also show
that many other standard constructions and vanishing theorems
produce $GV$-sheaves. A key technical point, already noted by Hacon
\cite{hacon}, is the very special nature of the Fourier-Mukai
transform of an ample line bundle on an abelian variety.

Let $X$ be smooth $d$-dimensional projective variety over a field of
characteristic $0$. Let $a: X \rightarrow A : = {\rm Alb}(X)$ be the
Albanese map of $X$, and assume that the dimension of $a(X)$ is
$d-k$ and the dimension of $A$ is $g$. Consider $\widehat A \cong
\Pic0$ and $\mathcal{P}$ a Poincar\' e line bundle on $A\times
\widehat A$. Consider also the pull-back $P :=(a\times {\rm
id}_{\widehat A})^*(\mathcal{P})$ on $X \times \widehat A$. In this
context one can define Mukai's Fourier functors for abelian
varieties
$$\R \hat{\mathcal{S}}:{\bf D}(A)\rightarrow {\bf D}(\widehat A) {\rm~ and~} \R \mathcal{S}:{\bf D}(\widehat A)\rightarrow {\bf D}(A)$$
given by the kernel $\mathcal{P}$ in both directions and, more
importantly for our purpose,
$${\bf R}\Phi_P :{\bf D}(X)\rightarrow {\bf D}(\widehat A) {\rm ~and~}
{\bf R}\Psi_P :{\bf D}(\widehat A)\rightarrow {\bf D}(X)$$ as in
\S2. Every $GV$-type condition in this section will be with respect
to $P$. Recall that for a nef $\QQ$-divisor $L$ on $X$, we denote by
$\kappa_L$ the Iitaka dimension $\kappa(L_{|F})$ of the restriction
to the generic fiber, if $\kappa(L) \ge 0$, and $0$ otherwise. We
are now ready to prove Theorem \ref{a} in the Introduction.

\begin{proof}(\emph{of Theorem \ref{a}.})
\noindent \emph{Step 1.} We first prove the Theorem in the case when
$D$ is an integral divisor, in other words we show the last
assertion. Let $L$ be a nef line bundle on $X$. It is enough to show
that $\omega_X\otimes L$ satisfies condition (2) in Theorem
\ref{general_f} (cf. also Corollary \ref{special}).

Let $M$ be ample line bundle on $\widehat A$, and assume for
simplicity that it is symmetric, i.e. $(-1_{\widehat A})^* M \cong
M$. We consider the two different Fourier transforms $\R \mathcal{S}
M=R^0S M$ (on A) and ${\bf R}\Psi_{P[g]} (M^{-1}) =R^g \Psi_P
(M^{-1}) =: \widehat{M^{-1}}$ (on X). These are both locally free
sheaves. We first claim that
$$ \widehat{M^{-1}} \cong  a^* R^g \SS (M^{-1}) \cong a^*(R^0 \mathcal{S} M)^\vee .$$
The second isomorphism follows from Serre duality, the symmetry of
$M$, and the fact that the Poincar\' e bundle satisfies the symmetry
relation $\mathcal{P}^\vee \cong ((-1_A)\times 1_{\widehat A})^*
\mathcal{P}$. If $a$ were flat, which is usually not the case, the
first isomorphism would follow simply by flat push-pull formula.
However, the same result holds since both $M$ and $R^gS (M^{-1})$
are locally free, hence flat. Less informally, consider the
cartesian diagram
$$\xymatrix{
X\times \widehat{A} \ar[r]^{a\times id} \ar[d]_{p_X} & A\times
\widehat{A}
\ar[d]^{p_A} \\
X \ar[r]^{a} & A } $$ Since $p_A$ is flat, one can apply e.g.
\cite{polishchuk1}, Theorem on p. 276, saying that there is a
natural isomorphism of functors
$${\bf L}a^*\circ {\bf R}{p_A}_* \cong {\bf
R}{p_X}_* \circ {\bf L}(a\times id_{\widehat A})^*.$$ The claim
follows since, as we are dealing with locally free sheaves, the two
derived pull-backs are the same as the usual
pull-backs.\footnote{Alternatively, in this particular case the
claim would follow in fact for any coherent sheaf, since the map on
the right is smooth -- cf. \cite{bo} Lemma 1.3.}

On the other hand, by \cite{mukai} 3.11, the vector bundle $R^0
\mathcal{S} M$ has the property:
$$\phi_M^*(R^0 \mathcal{S} M) \cong H^0(M)\otimes M^{-1}.$$
Here $\phi_M: \widehat A \rightarrow A$ is the standard isogeny
induced by $M$. We consider then the fiber product  $X^\prime
:=X\times_{A} \widehat A$ induced by $a$ and $\phi_M$:
$$\xymatrix{
X^\prime   \ar[r]^{\psi} \ar[d]_{b} & X
\ar[d]^{a} \\
\widehat{A} \ar[r]^{\phi_M} & A } $$ It follows that
\begin{equation}\label{square}
\psi^* \widehat{M^{-1}} \cong \psi^* a^* (R^0 \mathcal{S} M)^\vee
\cong b^*(H^0(M)\otimes M) \cong H^0(M)\otimes b^*M.
\end{equation}
Recall that we want to prove the vanishing of $H^i (\omega_X \otimes
L \otimes \widehat{M^{-1}})$ for $i > k - \kappa_L$. Since $\psi$,
like $\phi_M$, is \' etale, it is enough to prove this after
pull-back to $X^\prime$, so for $H^i (X^\prime, \omega_{X^\prime}
\otimes \psi^* L \otimes  \psi^* \widehat{M^{-1}})$. But by
$(\ref{square})$ we see that this amounts to the same for the groups
$H^i (X^\prime, \omega_{X^\prime} \otimes \psi^*L \otimes b^*{M})$.

We next use the fact that in this particular setting we have
additivity of Iitaka dimension in the following sense: 
$$\kappa(\psi^*L \otimes b^*M) \ge  d-k + \kappa_L.$$
In fact, assuming that $\kappa(L)\ge 0$, this is an actual equality
since $\psi^*L$ is nef and $b^*M$ is the pull-back of an ample line
bundle from a variety of dimension $d-k$, we have that
$$\kappa(\psi^*L \otimes b^*M) = \kappa (\psi^*L_{|F_b}) + d-k = \kappa(L_{|F}) + d-k,$$
where $F_b$ is the generic fiber of $b$ -- cf. \cite{mori},
Proposition 1.14.\footnote{The set-up for the definitions and
results in \emph{loc.cit.} is for fiber spaces over normal
varieties, but it is standard that the result can be easily reduced to 
that case by taking the Stein factorization.} On the other hand, if  
$\kappa(L) = -\infty$ it is still true that $\kappa(\psi^*L \otimes b^*M) \ge  d-k$. 
(By cutting with general hyperplane sections we can reduce to the case of generically
finite maps, when $\psi^*L \otimes b^*M$ is big and nef.)

Therefore the vanishing we want follows from the following
well-known extension (and immediate consequence) of the
Kawamata-Viehweg vanishing theorem for nef and big divisors (cf.
\cite{ev} 5.12):

\noindent $(*)$ \emph{Let $Z$ be a smooth projective variety and let
$N$ be a nef line bundle on $Z$, of Iitaka dimension $r\ge 0$. Then
$H^i(Z, \omega_Z\otimes N)=0$ for all $i>\dim Z - r$.}

\noindent \emph{Step 2.} Consider now, in the general case, $D$ an
effective $\QQ$-divisor on $X$ and $L$ a line bundle such that $L-D$
is nef. We keep the notation from the previous step. Let $f^\prime:
Y^\prime  \rightarrow X^\prime$ be a log-resolution of the pair
$(X^\prime , \psi^* D)$. Since $\psi$ is \' etale, by
\cite{positivity} 9.5.44 we have that $\J (\psi^* D)\cong \psi^{-1}
\J(D) \cong \psi^* \J (D)$, and as in Step 1 we are reduced to
showing the vanishing
$$H^i (X^\prime, \omega_{X^\prime} \otimes b^* M \otimes \psi^*L \otimes \J(\psi^* D) ) = 0,
{\rm ~for~all~} i> k.$$ The result follows from a
Nadel-vanishing-type version of $(*)$, which is the content of the
next Lemma.

\begin{lemma}\label{iitaka}
Let $X$ be a smooth projective complex variety of dimension $d$, $D$
an effective $\QQ$-divisor on $X$, and $L$ a line bundle such that
$L-D$ is nef and has Iitaka dimension $r\ge 0$. Then
$$H^i (X, \omega_{X} \otimes L \otimes \J(D) ) = 0, {\rm ~for~all~} i> d-r.$$
\end{lemma}
\begin{proof}
Like the proof of Nadel Vanishing, this is a standard reduction to
the integral Kawamata-Viehweg-type statement, and we only sketch it.
Take $f: Y \rightarrow X$ to be a log-resolution of the pair
$(X,D)$. By definition, $\J(D) = f_* \OO_Y( K_{Y/X} - [f^*D])$. By
Local Vanishing (cf. \cite{positivity} Theorem 9.4.1) and the
projection formula, using the Leray spectral sequence it is enough to 
have the desired vanishing for
$H^i (Y, \OO_Y(K_Y + f^*L - [f^*D]))$. As with the usual
Kawamata-Viehweg theorem (cf. \cite{positivity} 9.1.18), the
statement is reduced via cyclic covers to the integral case, which
is then covered by $(*)$.
\end{proof}
This concludes the proof of the Theorem.
\end{proof}

For applications it is useful to note that the statement of Theorem
\ref{a} can be extended to the setting of asymptotic multiplier
ideals. In this case we do not need the nefness assumption. For a
line bundle $L$ with $\kappa(L)\ge 0$, we denote by $\J(\parallel L
\parallel)$ the asymptotic multiplier ideal associated to $L$ (cf.
\cite{positivity} \S 11.1).

\begin{corollary}\label{asymptotic}
Let $X$ be a smooth projective variety of dimension $d$ and Albanese
dimension $d-k$, and $L$ be a line bundle on $X$ with
$\kappa(L)\ge 0$. Then $\omega_X \otimes L\otimes \J(\parallel
L\parallel)$ is a $GV_{-(k-\kappa_L)}$-sheaf.
\end{corollary}
\begin{proof}
This is a corollary of the proof of Theorem \ref{a}. By the behavior
of asymptotic multiplier ideals under \'etale covers (cf.
\cite{positivity}, Theorem 11.2.23), we have that $\psi^*
\J(\parallel L \parallel) \cong \J(\parallel \psi^* L \parallel)$.
As before, we need to show that for a sufficiently positive line
bundle $M$ on $A$ we have
$$H^i (X^{\prime}, \omega_{X^{\prime}} \otimes \psi^* L \otimes b^*M \otimes
\J (\parallel \psi^* L \parallel) ) = 0, ~\forall ~i > k -
\kappa_L.$$ By definition we have $\J(\parallel \psi^*L \parallel) =
\J(\frac{1}{p} \cdot |p\psi^*L|)$, where $p$ is a sufficiently large
integer. We consider $f: Y \rightarrow X^{\prime}$ a log-resolution
of the base locus of the linear series $|p\psi^*L|$, and write
$f^*|p\psi^*L| = |M_p| + F_p$, where $M_p$ is the moving part (in
particular nef), and $F_p$ is the fixed part. Thus $\J(\parallel
\psi^*L \parallel) = f_*\OO_Y ( K_{Y/X^{\prime}} - [\frac{1}{p}
F_p]),$ so by Local Vanishing and the projection formula, using the Leray
spectral sequence we have
$$H^i (X^{\prime}, \omega_{X^{\prime}} \otimes \psi^* L \otimes b^*M \otimes
\J (\parallel \psi^*L \parallel) )\cong H^i (Y , \OO_Y( K_Y +
f^*\psi^*L - [\frac{1}{p} F_p] + f^*b^*M)).$$ Note now that
$f^*\psi^* L - [\frac{1}{p} F_p]$ is numerically equivalent to the
$\QQ$-divisor $\frac{1}{p} M_p + \{\frac{1}{p} F_p\}$, where the
last factor is fractional and simple normal crossings, so with
trivial multiplier ideal. Since $M_p$ is nef, as in Theorem \ref{a}
we have that $H^i = 0$ for $i > d - \kappa (M_p + f^*b^*M)$. On the
other hand $M_p$ asymptotically detects all the sections of
$\psi^*L$, and the generic fiber for $b$ and $b\circ f$ is the same,
so we have that
$$\kappa (M_p + f^*b^*M) = d-k  + \kappa_{M_p} = d-k + \kappa_L.$$
\end{proof}

\begin{corollary}\label{canonical}
Let $X$ be a smooth projective variety whose Albanese map has
generic fiber of general type (including the case of maximal
Albanese dimension). Then $\OO((m+1)K_X) \otimes \J(\parallel mK_X
\parallel)$ is a $GV$-sheaf for all $m\ge 1$.
\end{corollary}

\begin{remark}
(1) Theorem \ref{a} holds more generally, but with the same proof,
replacing the Albanese map with any morphism to an abelian variety
$a:X\rightarrow A$.

\noindent (2) A particular case of Theorem \ref{a}, extending the
theorem of Green-Lazarsfeld to the case of line bundles of the form
$\omega_X \otimes L$ with $L$ semiample, or even nef and abundant,
was proved by H. Dunio. This was done by reducing to the case of the
Green-Lazarsfeld theorem via cyclic covers (cf. \cite{ev} 13.7 and
13.10). Note also that Ch. Mourougane (cf. \cite{mourougane},
Th\'eorem\`e) showed using similar covering techniques that the same
is true on a compact K\" ahler manifold. The main restriction in
both cases essentially has to do, at least asymptotically, with the
existence of lots of sections. We expect that Theorem \ref{a} is
also true on compact K\" ahler manifolds $(X,\omega)$ if $L$ is only
assumed to be nef, which in this case means that $c_1(L)$ is in the
closure of the K\"ahler cone, i.e. the closed cone generated by
smooth non-negative closed $(1,1)$-forms.

\noindent (3) Results on the structure of cohomological support loci
for twists with multiplier ideal sheaves, under some specific
numerical hypotheses, are contained in \cite{budur}.
\end{remark}

\noindent {\bf Pluricanonical bundles.} Generic Vanishing for
$\omega_X$ states that if the Albanese dimension of $X$ is $d-k$,
then $\omega_X$ is a $GV_{-k}$-sheaf.  Theorem \ref{a} implies that the
same is true for all powers $\omega_X^{\otimes m}$ when $X$ is
minimal. In fact, as soon as $m\ge 2$, one can do better according
to the Kodaira dimension of the generic fiber of $a$. Recall that we
denote by $\kappa_F$ the Kodaira dimension $\kappa(F)$ if this is
non-negative, or $0$ if $\kappa(F) = -\infty$.

\begin{corollary}
Let $X$ be a smooth projective minimal variety. Then for any nef
line bundle $L$ on $X$ and any $m\ge 1$, $\omega_X^m \otimes L$ is a
$GV_{-(k- \kappa_F)}$-sheaf. In particular $\omega_X^{\otimes m}$ is
a $GV_{-(k- \kappa_F)}$-sheaf for all $m \ge 2$.
\end{corollary}
\begin{proof}
The minimality condition means that $\omega_X$ is nef. Everything
follows directly from the proof Theorem \ref{a}, together with the
extra claim that $\kappa(F) = -\infty$ if and only if $\kappa(K_X +
b^*M) = -\infty$. It is well known (cf. \cite{mori}, Proposition
1.6), that $\kappa(F) = -\infty$ implies $\kappa(K_X + b^*M) =
-\infty$. The statement that $\kappa(F) \ge 0$ implies $\kappa(K_X +
b^*M) \ge 0$ follows from the general results on the subadditivity
of Kodaira dimension (cf. \emph{loc.cit.}).
\end{proof}

As J. Koll\'ar points out, it is very easy to see that if the
minimality condition is dropped, then the higher powers
$\omega_X^{\otimes m}$, $m\ge 2$, do not satisy generic vanishing.

\begin{example}\label{surfaces}
Let $Y$ be the smooth projective minimal surface of maximal Albanese
dimension (for example an abelian surface), and $f: X \rightarrow Y$
the blow-up of $Y$ in a point, with exceptional divisor $E$.  By the
previous Corollary $\omega_Y^{\otimes m}$ is a $GV$-sheaf for all
$m$. We claim that this cannot be true for $\omega_X^{\otimes m}$
when $m \ge 2$. Indeed, if one such were $GV$, then we would have
that
$$\chi (\omega_X^{\otimes m}) = h^0 (\omega_X^{\otimes m} \otimes f^* \alpha) {\rm ~~and~~} \chi (\omega_Y^{\otimes m}) = h^0 (\omega_Y^{\otimes m} \otimes \alpha)$$
for $\alpha \in {\rm Pic}^0(Y)$ general. The fact that $\omega_X
\cong f^* \omega_Y \otimes \OO_Y(E)$ implies easily that $h^0
(\omega_X^{\otimes m} \otimes f^* \alpha)= h^0 (\omega_Y^{\otimes m}
\otimes \alpha)$. On the other hand, the Riemann-Roch formula shows
that $\chi (\omega_X^{\otimes m}) \neq \chi (\omega_Y^{\otimes m})$
as soon as $m \ge 2$, which gives a contradiction.
\end{example}

\begin{remark}
The previous example shows that in general the tensor product of
$GV$-sheaves is not $GV$. This is however true on abelian varieties,
cf. \cite{pp3}.
\end{remark}

\noindent {\bf Higher direct images.} 
In analogy with work of Hacon (see
\cite{hacon} and also \cite{hp} -- cf. Remark
\ref{previous_higher} below), we show that higher direct images of
dualizing sheaves have the same generic vanishing behavior as
dualizing sheaves themselves. In order to have
  the most general statement, we replace the Albanese map of a smooth variety with the following setting: we consider $X$ to be an arbitrary Cohen-macaulay projective variety, and $a: X
\rightarrow A$ a morphism to an abelian variety. We discuss the
$GV$-property with respect to the Fourier-Mukai functor induced by
the kernel $P = (a\times id)^* \mathcal{P}$ on $X \times \widehat
A$, where $\mathcal{P}$ is the Poincar\'e bundle of $A \times
\widehat A$.

\begin{theorem}\label{direct}
Let $f:Y \rightarrow X$ be a morphism, with $X$, $Y$
projective, $Y$ smooth and $X$ Cohen-Macaulay. Let $L$ be a nef line bundle on $f(Y)$ (reduced image of $f$). If
the dimension of $f(Y)$ is $d$ and that of $a(f(Y))$ is $d-k$, then $R^j
f_* \omega_Y \otimes L$ is a $GV_{-(k - \kappa_L)}$-sheaf on $X$ for
any $j$.
\end{theorem}
\begin{proof}
This follows along the lines of the proof of Theorem \ref{a}, so we
only sketch the argument. Consider  a sufficiently ample line bundle
$M$ on $\widehat A$, and form the two cartesian squares:
$$\xymatrix{
Y^\prime \ar[d]^{\nu} \ar[r]^{f^\prime} & X^\prime   \ar[r]^{b} \ar[d]^{\psi} & \widehat A \ar[d]^{\phi_M} \\
Y \ar[r]^{f} & X \ar[r]^{a} & A } $$ where $\phi_M$ is the standard
isogeny induced by $M$. As before, we need to check the vanishing
$$H^i (X^\prime, R^jf^{\prime}_* \omega_{Y^\prime} \otimes \psi^*L \otimes b^*M) = 0,
{\rm ~for~all~} i >k {\rm~and~all~} j.$$ We have that $\psi^*L
\otimes b^*M$ is nef and has $\kappa (\psi^*L \otimes b^*M) \ge
\kappa_L + d-k$. The required vanishing is a consequence of the
variant of Koll\' ar's vanishing theorem (cf. \cite{kollar}, Theorem
2.1(iii)) stated in the Lemma below. (Note that the Lemma is applied replacing $X^\prime$ with  $f^\prime(Y^\prime)$).
\end{proof}

\begin{lemma}\label{higher_direct}
Let $f:Y \rightarrow X$ be a  surjective morphism of projective
varieties, with $Y$ smooth and $X$ of dimension $d$. If $M$ is a nef
line bundle on $X$, of Kodaira-Iitaka dimension $\kappa (M) = d-k$,
then
$$H^i (X, R^j f_* \omega_Y \otimes M) = 0, {\rm ~for~all~} i>k {\rm ~and ~all~} j.$$
\end{lemma}
\begin{proof}
This uses the full package provided by \cite{kollar} Theorem 2.1. We
note to begin with that the condition on $M$ implies that there
exist hypersurfaces $H_1, \ldots, H_k$ on $X$ such that $Z : =  H_1
\cap \ldots \cap H_k \subset X$ is a subvariety of dimension $d-k$
such that $M_{|Z}$ is big and nef. We can assume that the $H_i$'s
are sufficiently positive and general, so in particular by Bertini 
the preimages $\tilde H_i$ of the $H_i$'s in $Y$ are smooth. 

We do a descending induction: let $H$ in $X$ be one of the hypersurfaces $H_i$ as above. Pushing
forward the obvious adjoint sequence on $Y$, we obtain a long exact sequence:
$$\ldots \rightarrow R^{j-1}f_*  \omega_{\tilde H} \rightarrow R^{j}f_*  \omega_Y \rightarrow
R^{j}f_*  (\omega_Y (\tilde H))\rightarrow R^{j}f_*
\omega_{\tilde H}\rightarrow R^{j+1}f_*  \omega_Y
\rightarrow \ldots$$ 
The sheaves $R^{j}f_*  \omega_{\tilde H}$ are
supported on $H$, while by  \cite{kollar} Theorem 2.1(i), the
sheaves $R^{j}f_*  \omega_Y$ are torsion-free. This implies
that the long exact sequence above breaks in fact into short exact
sequences
$$0 \rightarrow R^{j}f_*  \omega_Y\rightarrow
R^{j}f_*  (\omega_Y (\tilde H))\rightarrow  R^{j}f_*
\omega_{\tilde H}\rightarrow 0.$$ 
We twist these exact sequences by
$M$, and pass to cohomology. Since $H$ can be taken sufficiently
positive and $f_* \OO_Y (\tilde H) \cong \OO_X(H)$, by Serre vanishing we may assume that $H^i (R^{j}f_*
(\omega_{\tilde Y} (\tilde H)) \otimes M) = 0$ for all $i >0$. This
implies that $H^i (R^{j}f_* \omega_Y \otimes M) \cong
H^{i-1} (R^{j}f_* \omega_{\tilde H} \otimes M)$ for all $i$ and all
$j$. We continue intersecting with $H_i$'s until we get to $Z$. This
implies that
$$H^i (X, R^{j}f_* \omega_Y \otimes M) \cong H^{i-k} (Z, R^{j}f_* \omega_{\tilde Z} \otimes M_{|Z}).$$
But this is $0$ for $i >k$, since on $Z$ we can apply the vanishing
theorem \cite{kollar} Theorem 2.1(iii) (or rather its well-known
version for big and nef line bundles).
\end{proof}

\begin{remark}\label{previous_higher}
In the case when one considers higher direct images $R^j f_*
\omega_Y$, i.e. $L = \OO_X$, Theorem \ref{direct}  was already noted
in \cite{hacon} Corollary 4.2, and more generally \cite{hp} Theorem
2.2(a) and the references therein.
\end{remark}

\noindent {\bf Generic Nakano-type vanishing.} Using similar
techniques, we can deduce generic vanishing results for bundles of
holomorphic forms, based on a suitable generalization of Nakano
vanishing.

\begin{theorem}\label{nakano_1}
Let $X$ be a smooth projective variety, with Albanese image of
dimension $d-k$. Denote by $m$ the \emph{maximal} dimension of a
fiber of $a$, and consider $l:= {\rm max}\{k, m-1\}$. Then:\\
(1) $\Omega_X^j$ is a $GV_{-(d-j+l)}$-sheaf for all $j$.\\
(2)  ${\rm codim}_{\hat A} V^i (\Omega_X^j) \ge {\rm max}\{i + j - d
-l, d - i - j -l\}$, for all
 $i$ and all $j$.\\
(3) If $L$ is a nef line bundle on $X$ and $a$ is \emph{finite},
then $\Omega^j_X\otimes L$ is a $GV_{-(d - j)}$-sheaf for  all $j$.
\end{theorem}
\begin{proof}
 (1) Again we follow precisely the pattern of the proof of Theorem
\ref{a}, replacing $\omega_X \otimes L$ with $\Omega_X^j$. The
result reduces to checking the vanishing
$$H^i (X^\prime, \Omega_{X^\prime}^j \otimes b^*M) = 0, ~{\rm for ~all~}
i > d - j + l.$$ But, as the pull-back of an ample line bundle via
$b$, the line bundle $b^*M$ is $m$-ample in the sense of Sommese (\cite{ev} 6.5).
Thus the needed vanishing follows from the generalization by Sommese
and Esnault-Viehweg of the Nakano vanishing theorem: if $L$ is an
$m$-ample line bundle on $Y$ of dimension $d$, then $H^i (Y, \Omega_Y^j
\otimes L) = 0$ for all $i > d - j  + {\rm max}\{d - \kappa(L), m -
1\}$ (cf. \cite{ev} 6.6).\\
(2) Apply part (1) to the sheaves $\Omega_X^j$ and
$\Omega_X^{d-j}$, which are related in an obvious way by Serre
duality.\\
(3)
 As before, we are reduced to checking the vanishing
of the cohomology groups $H^i (X^\prime, \Omega_{X^\prime}^j \otimes
\psi^*L \otimes b^*M)$. But since $b$ is a finite map, $\psi^*L
\otimes b^*M$ is ample, and so the result follows from the Nakano
vanishing theorem.
\end{proof}

\begin{remark}
Note that Nakano vanishing does not hold in the more general setting
of twisting with big and nef line bundles (cf. \cite{positivity},
Example 4.3.4). Thus one does not expect to have generic vanishing
for $\Omega_X^j$ depending only on the dimension of the generic
fiber of the Albanese map, as in the case of $\omega_X$. A
counterexample was indeed given by Green and Lazarsfeld (\cite{gl1},
Remark after Theorem 3.1). The above shows that there are however
uniform bounds depending on the \emph{maximal} dimension of the
fibers of the Albanese map -- in particular if the Albanese map is
equidimensional, or if the fiber dimension jumps only by one, then
the exact analogue of generic vanishing for $\omega_X$ does
hold.\footnote{For example, in the case of finite Albanese maps, the
results of \cite{gl1} seem to imply only the weak form of generic
vanishing, as in our Corollary \ref{b} in the Introduction.} Note
that in the counterexample mentioned above, the fiber dimension
jumps by $2$: it is shown there that $\Omega_X^1$ does not satisfy
the expected $GV$ condition. The result above shows that it does
satisfy the ``one worse" $GV$ condition. Finally, in \cite{gl1}
Theorem 2, another variant for generic Nakano-type vanishing is
proposed in terms of zero-loci of holomorphic $1$-forms. It would 
be interesting if the two approaches could be combined.
\end{remark}

\noindent {\bf Vector bundles.} By the same token, the various known
vanishing theorems for higher rank vector bundles show that nef
vector bundles also satisfy a weaker form of generic vanishing.
Below is a sampling of results. All the notions and results we refer to can be
found in \cite{positivity} \S7.3.

\begin{theorem}\label{vector}
Let $X$ be a smooth projective variety of dimension $d$ and $E$ a
nef vector bundle on $X$. Then:
\newline
\noindent (1) If the Albanese map of $X$ is \emph{finite}, then \
$\omega_X\otimes \Lambda^a E$ is a $GV_{-({\rm rk}(E)-a)}$-sheaf.
Moreover \ $\Omega_X^j\otimes E$ is a $GV_{-({\rm
rk}(E)+d-j-1)}$-sheaf for all $j$.\\
 (2) If the Albanese map of $X$ is
\emph{generically finite}, then for all $m \ge 0$, \ $\omega_X\otimes
S^m E \otimes {\rm det}(E)$ is a $GV$-sheaf. Moreover, if $E$ is
$k$-ample in the sense of Sommese, then \ $\omega_X\otimes \Lambda^a
E$ is a $GV_{-({\rm rk}(E)+k-a)}$-sheaf for all $a>0$.
\end{theorem}
\begin{proof}
Everything goes exactly as in the proof of Theorem \ref{a}, so in
the end one is reduced to checking vanishing for cohomology groups
of the form $H^i (\omega_X \otimes F\otimes L)$ and the
corresponding Akizuki-Nakano analogues, where $F$ is a vector bundle
as above, and $L$ is an ample or big and nef line bundle. Then (1)
follows from the Le Potier vanishing theorem, the second part of (2)
follows from Sommese's version of the same theorem, while the first
part of (2) from the Griffiths vanishing theorem (which also works
in the big and nef case -- cf. \cite{positivity} 7.3.2 and 7.3.3).
\end{proof}

Further, more refined results along these lines, and especially
taking into account the beautiful general vanishing theorems for
vector bundles of Demailly, Manivel, Arapura and others, can be
formulated  by the interested reader. In the K\"ahler case, results
on Nakano semipositive vector bundles with some stronger conditions
on the twists were proved by Mourougane \cite{mourougane}.

\noindent {\bf Algebraicity and positive characteristic.} The
methods of this paper give an algebraic proof of Generic Vanishing
in characteristic $0$. The main point is that the Kawamata-Viehweg
theorem can be reduced in characteristic $0$ to Kodaira vanishing,
via covering constructions. This in turn is proved via reduction mod
$p$ in \cite{di}. We will comment later that other known results can
be similarly proved algebraically -- cf. Remark \ref{algebraic}.

On the other hand, assume that $X$ is defined over a perfect field
of characteristic ${\rm char }(k) \ge {\rm dim}(X)$, and that it
admits a lifting to the $2$nd Witt vectors $W_2(k)$. Then, again by
\cite{di}, Kodaira vanishing is still known to hold. However, this
is not (yet) the case with the analogue of Kawamata-Viehweg (cf.
\cite{ev} 11.6 and 11.7) -- it would follow as in characteristic $0$
if we had embedded resolution of singularities over the field $k$
and over $W_2(k)$.

As a consequence, with the current state of knowledge we know that
the main results on Generic Vanishing in this paper (especially
Theorem \ref{a}) hold in positive characteristic, under the
assumptions above, only if either one of the following holds:\\
 (1) The Albanese map is finite onto its image. \hfill\break
 (2) The dimension of $X$ is at most three (cf. \cite{abhyankar}),
if ${\rm char}(k) > 5$ and the embedded resolution also admits a
lifting over $W_2(k)$.\hfill\break
 (3) The standard generic
vanishing for $\omega_X$ holds in arbitrary characteristic if the
Albanese map is separable, by a result of the first author
\cite{pareschi}.

\section{Applications via the Albanese map}

We give some examples of how Theorems \ref{a} and \ref{direct} can
be applied to basic questions in the spirit of \cite{kollar2} and
\cite{el}, and we make another comment on algebraic proofs.

\noindent {\bf Existence of sections and generic plurigenera.} The
first is related to the existence of sections. Recall that it is
known that every nef line bundle on an abelian variety is
numerically equivalent to an effective one. Also, if $L$ is a nef
and big line bundle on a variety of maximal Albanese dimension, then
$K_X + L$ has a non-zero section.\footnote{Cf. e.g. \cite{pp2} \S5
for a quick proof of this; in fact much more holds, cf.
\cite{kollar} Theorem 16.2.}

\begin{theorem}\label{sections}
Let $X$ be a smooth projective variety of maximal Albanese
dimension. Let $L$ be a line bundle on $X$ such that either one of
the following holds:

\indent (1) $L$ is nef.\hfill\break \indent (2) $\kappa(L) \ge 0$.

\noindent Then there exists $\alpha\in \Pic0$ such that
$h^0(\omega_X \otimes L\otimes \alpha) > 0$. In particular $K_X + L$
is numerically equivalent to an effective divisor.
\end{theorem}
\begin{proof}
Assume first that $L$ is nef. If the conclusion doesn't hold, we
have that $V^0 (\omega_X \otimes L) = \emptyset$. But since by the
Theorem $\omega_X \otimes L$ is $GV$, Proposition
\ref{general_inclusions} implies that $V^i(\omega_X \otimes L) =
\emptyset$  for all $i$. By Grauert-Riemenschneider vanishing we
have that $R^j a_* (\omega_X\otimes L) = 0$ for all $j >0$, so $V^i
(a_* (\omega_X \otimes L)) = V^i (\omega_X \otimes L) $ for all $i$.
This means that $a_* (\omega_X \otimes L)$ is a non-zero (since $a$
is generically finite) sheaf on $A$ whose Fourier-Mukai transform
$\R \hat{\mathcal{S}}(a_*(\omega_X \otimes L))$ is equal to zero.
But this is impossible since $\R \hat{\mathcal{S}}$ is an
equivalance.

If $L$ has non-negative Iitaka dimension, we can consider the
asymptotic multiplier ideal $\J(\parallel L \parallel)$ as in
Corollary \ref{asymptotic}, and we show that in fact there exists
$\alpha\in \Pic0$ such that
$$h^0(\omega_X \otimes L\otimes \J(\parallel L \parallel) \otimes \alpha) > 0.$$
As above, by Corollary \ref{asymptotic} $\omega_X \otimes L\otimes
\J(\parallel L \parallel)$ is a $GV$-sheaf, so assuming that the
conclusion doesn't hold, we get a contradiction if we know that $V^i
(a_*(\omega_X \otimes L\otimes \J(\parallel L \parallel))  ) = V^i
(\omega_X \otimes L\otimes \J(\parallel L \parallel))$ for all $i$,
which in turn follows if we know that
$$R^j a_* (\omega_X \otimes L\otimes \J(\parallel L \parallel))   = 0 {\rm~for~all~} j >0. $$
Recall that $\J(\parallel L \parallel) = \J(\frac{1}{m}\cdot |mL|)$
for some $m \gg0$, and consider $\phi: Y \rightarrow X$ a
log-resolution of the base locus of $|mL|$. Write $\phi^* (mL) = M_m
+ F_m$, where $M_m$ is the moving part and $F_m$ the fixed part, in
simple normal crossings. We see easily that
$$R^j a_* (\OO_X( K_X+ L)\otimes \J(\parallel L \parallel))  \cong
 R^j (a\circ \phi)_* \OO_Y ( K_Y + N),$$
 where $N \equiv_{\QQ} \frac{1}{m} M_m + \{\frac{1}{m} F_m\}$. Since $M_m$ is nef and $a\circ \phi$ is
 generically finite,  this follows again from the ($\QQ$-version of) Grauert-Riemenschneider-type vanishing 
 (cf. \cite{km} Corollary 2.68, noting that the same proof works for generically finite maps).
\end{proof}

Next we use Corollary \ref{canonical} for a result which
interpolates between Koll\' ar's theorem on the multiplicativity of
plurigenera under \'etale maps for varieties of general type
(cf. \cite{kollar2} Theorem 15.4), and the (obvious) case of abelian
varieties. The invariant which is well-behaved under \'etale covers
in this case is the \emph{generic plurigenus}:
$$P_{m, gen} : = h^0 (X, \OO(mK_X)\otimes \alpha),$$
where $\alpha\in \Pic0$ is taken general enough so that the quantity
is minimal over $\Pic0$.

\begin{theorem}\label{plurigenera}
Let $Y \rightarrow X$ be an \'etale map of degree $e$ between smooth
projective varieties whose Albanese maps have generic fiber of
general type (including the case of maximal Albanese dimension).
Then
$$P_{m, gen}(Y) =  e\cdot P_{m, gen}(X), ~{\rm for~all} ~m\ge2.$$
\end{theorem}
\begin{proof}
The proof follows Koll\'ar idea of expressing the plurigenus as an
Euler characteristic, but in the language of asymptotic multiplier
ideals as in \cite{positivity} Theorem 11.2.23. Fix $m \ge 2$ and
$\alpha \in \Pic0$ sufficiently general so that
$$P_{m, gen}(X) = h^0 (X, \OO(mK_X)\otimes \alpha).$$
Since torsion points are dense in $\Pic0$, we can further assume
that $\alpha$ is  torsion. The point is that asymptotic multiplier
ideals do not detect torsion. Indeed:
$$\J( \parallel mK_X + \alpha \parallel ) \cong \J(\parallel mK_X \parallel),$$
since these are  computed from the linear series $|p(mK_X +
\alpha)|$ and $|pmK_X|$, for any $p$ sufficiently large, and in
particular divisible enough so that it kills $\alpha$.   On the
other hand, we know that
$$H^0(\OO(mK_X) \otimes \alpha )\cong H^0(\OO(mK_X) \otimes \alpha \otimes
\J( \parallel mK_X + \alpha \parallel )), $$ and also that
$\J(\parallel mK_X \parallel) \subseteq \J(\parallel (m-1)K_X
\parallel)$, so as a consequence we have
$$P_{m, gen}(X) = h^0 (\OO(mK_X) \otimes \alpha \otimes \J(\parallel (m-1)K_X \parallel)).$$
At this stage we can use Generic Vanishing: by Corollary
\ref{canonical}, we know that the sheaf $\OO(mK_X) \otimes
\J(\parallel (m-1)K_X \parallel)$ is $GV$, so $\alpha$ could also be
chosen such that
$$H^i (\OO(mK_X) \otimes \alpha \otimes \J(\parallel (m-1)K_X \parallel)) = 0, ~\forall~i>0.$$
This finally implies that
$$P_{m, gen}(X) = \chi (\OO(mK_X) \otimes \J(\parallel (m-1)K_X \parallel)),$$
since Euler characteristic is invariant under deformation. Since the
same is true for $Y$, the result follows immediately from the
multiplicativity of Euler characteristics under \'etale maps, and
the fact that
$$ \J(\parallel (m-1)K_Y \parallel) \cong f^* \J(\parallel (m-1)K_X \parallel),$$
which is a consequence of the behavior of asymptotic multiplier
ideals under \'etale covers, \cite{positivity} Theorem 11.2.16.
\end{proof}

This implies subadditivity of generic plurigenera, precisely as in
\cite{kollar} Theorem 15.6.

\begin{corollary}\label{subadditivity}
Let $X$ be  a smooth projective variety with nontrivial algebraic
fundamental group, whose Albanese map has generic fiber of general
type (for example a variety of maximal Albanese dimension). Then,
for $m, n\ge 2$, if $P_{m, gen}(X) > 0$ and $P_{n, gen}(X) > 0$,
then
$$P_{m+n, gen}(X) \ge P_{m, gen}(X)  + P_{n, gen}(X).$$
\end{corollary}

Since a variety $X$ of maximal Albanese dimension has generically
large fundamental group (see \cite{kollar2} Definition 4.6 for the slightly technical definition), 
by \cite{kollar2} Theorem 16.3 we have that
if $X$ is of general type, then $P_m (X) \ge 1$ for $m \ge 2$ and
$P_m (X) \ge 2$ for $m \ge 4$. Using the above, a stronger statement
can be made.

\begin{corollary}
Let $X$ be a variety of maximal Albanese dimension and of general
type. Then
\newline
\noindent (1) $P_{m, gen}(X) \ge 1$ for $m \ge 2$.
\newline
\noindent (2) $P_{m, gen}(X) \ge 2$ for $m \ge 4$.
\end{corollary}
\begin{proof}
The first statement is already known, in a more general form.
Indeed, since bigness is preserved under numerical equivalence, we
have that $\omega_X \otimes \alpha$ is big for all $\alpha \in
\Pic0$. But on varieties of maximal Albanese dimension (or more
generally with generically large fundamental group), every line
bundle of the form $\omega_X \otimes L$ with $L$ big has non-zero
sections (cf. \cite{kollar2} Theorem 16.2). The second follows from
(1) and Corollary \ref{subadditivity}.
\end{proof}

\noindent {\bf Components of $V^0$.} Let $X$ be a smooth projective
variety of dimension $d$, and $a:X \rightarrow A$ its Albanese map.
We denote as always the dimension of the generic fiber of $a$ by
$k$. We can assume without loss of generality that the image $Y =
a(X)$ is smooth by passing to a resolution of singularities, which
is sufficient to check the dimension properties we are interested
in. First of all, Corollary \ref{isolated} gives in this case:

\begin{corollary}\label{surjective}
Say $X$ is of maximal Albanese dimension, and there exists a 
$GV$-sheaf $\F$ on $X$ with an isolated point in $V^0(\F)$.
Then ${\rm dim}~X \ge {\rm dim}~A$, and so the Albanese map is
surjective and $\dim X = \dim A$.
\end{corollary}

A special case is the following extension of \cite{el} Proposition
2.2; cf. also Remark \ref{algebraic}.

\begin{corollary}\label{old}
Assume that the characteristic is zero, and that there exists an
isolated point in $V^0 (\omega_X\otimes a^*L)$ for some nef line
bundle $L$ on $a(X)$. Then the Albanese map $a$ is surjective.
\end{corollary}
\begin{proof}
By passing to a resolution of singularities, we can assume that the 
Albanese image $Y = a(X)$ is smooth. 
Theorem \ref{direct} implies that $a_* \omega_X\otimes L$ is a
$GV$-sheaf on $Y$. By hypothesis there is an isolated
point in $V^0 (a_* \omega_X \otimes L)$, so Corollary
\ref{surjective} applies to give $Y = A$.
\end{proof}

The main result is that Theorem \ref{direct} and Proposition
\ref{component}, together with the structure result for $V^0
(\omega_X)$ in \cite{gl2}, imply the generalization of \cite{el}
Theorem 3 mentioned in the Introduction.

\begin{proof}(\emph{of Theorem \ref{e}}.)
Assume that there is a component $W$ of $V^0 (\omega_X)$ of
codimension $p > 0$. By \cite{gl2} we know that $W$ is a translation
of an abelian subvariety of $\widehat A = \Pic0$, which we will abusively also denote
by $W$. Denote $C: = \widehat W$, the dual abelian variety, and
consider the sequence of homomorphisms of abelian varieties
$$1 \rightarrow B \rightarrow A \rightarrow C \rightarrow 1.$$

Let $Y = a(X)$, and consider the morphism $f: Y \rightarrow C$
induced by the composition of the inclusion in $A$ and the
projection to $C$. Denote by $k$ the dimension of the generic fiber
of $f$. To prove the Theorem it is enough to show that $k \ge p$.
Indeed, the fibers of $A \rightarrow C$ are subtori of $A$ of
dimension $p$.

Now Theorem \ref{direct} implies that $a_* \omega_X$ is a $GV$-sheaf
on $Y = a(X)$, and we clearly have $V^0 (a_* \omega_X) = V^0
(\omega_X)$. This, together with Proposition \ref{component},
implies that $W$ is also a component of $V^p (a_* \omega_X)$. On the
other hand, Theorem \ref{direct} also implies that $a_* \omega_X$ is
a $GV_{-k}$-sheaf with respect to the natural Fourier-Mukai transform
$\D(Y) \rightarrow \D(C)$, so that
$${\rm codim}~V^p_{\widehat C} (a_* \omega_X) \ge p -k.$$
But $\widehat C =  W$, so by definition the line bundles in $\Pic0$
parametrized by $W$ are precisely those pulled back from $C$. We
finally have that $W \subseteq V^p_{\widehat C} (a_* \omega_X)
\subseteq \widehat C = W$, which implies equality everywhere. Hence
the codimension of $V^p_{\widehat C} (a_* \omega_X)$ is in fact $0$,
which gives $k \ge p$.
\end{proof}

\begin{remark}[Algebraic proofs]\label{algebraic}
We note that in the present approach one does not need to appeal to
complex analytic techniques. For instance, the Theorem above gives
in particular an algebraic proof of \cite{el} Theorem 3. Moreover,
the case $L = \OO_X$ in Corollary \ref{old}, together with the
results of Pink-Roessler \cite{pr}, provides an algebraic proof of
another result of Ein-Lazarsfeld saying that if $P_1(X) = P_2(X)
=1$, then the Albanese map of $X$ is surjective (cf. \cite{el}
Theorem 4). This in turn implies the same for Kawamata's well-known
theorem \cite{kawamata} saying that if $\kappa (X) = 0$, then the
Albanese map is surjective. Indeed, in \cite{pr} it is shown via the
reduction mod $p$ method of Deligne-Illusie, that the $V^j
(\omega_X)$ are unions of translates of subtori of $\Pic0$. This
allows the proof of \cite{el} Proposition 2.1 to go through
unchanged, while Corollary \ref{old} also shows that \cite{el}
Proposition 2.2 has an algebraic proof.
\end{remark}

\section{Applications and examples for bundles on curves and Calabi-Yau fibrations}

One of the main features of Theorem \ref{general_f} is that it
applies to essentially any integral transform. Here we exemplify
with some statements for vector bundles on curves and on some
threefold Calabi-Yau fibrations.

\noindent {\bf Semistable vector bundles on curves having a theta
divisor.} Let $X$ be a smooth projective curve of genus $g \ge 2$,
and let $SU_X(r, L)$ be the moduli space of semistable vector
bundles on $X$ of rank $r$ and fixed determinant $L \in {\rm
Pic}^d(X)$. The Picard group of $SU_X(r, L)$ is generated by the
determinant line bundle $\LL$. Results in this paper and the Strange
Duality provide a Fourier-Mukai criterion for detecting base points
of the linear series $|\LL^k|$ for all $k$.

We start with a special case when this base locus is much better
understood, namely the case of $d = r(g-1)$ and $k=1$. In this case
at least part of the criterion was already noted by Hein \cite{hein}
(see below). It is well known (cf. \cite{beauville} \S3) that if $d
= r(g-1)$, then a semistable vector bundle $E$ (or rather its
$S$-equivalence class) is in the base locus of $|\LL|$ if and only
if
$$H^0 (E \otimes \xi) \neq 0, {\rm ~for~all~} \xi \in \Pic0.$$
Otherwise, the locus described set-theoretically as $\Theta_E := \{
\xi ~|~ h^0 (E \otimes \xi) \neq 0\} \subset \Pic0$ is a divisor,
and one says that $E$ \emph{has a theta divisor}. We have the
elementary:

\begin{lemma}\label{base1}
If $E$ is as above, then the following are equivalent:\\
(1) $E$ is not a base point for $|\LL|$ (i.e. $E$ has a theta
divisor).\\
(2) $E$ is a $GV$-sheaf with respect to any Poincar\'e bundle $P$ on
$X \times \Pic0$.\\
(3) $R^0 \Phi_P E = 0$, i.e. $E$ satisfies $WIT_1$ with respect to
$P$.
\end{lemma}
\begin{proof}
The equivalence of (1) and (2) is the above discussion, plus the
fact that since $\chi(E) = 0$,  for all $\xi \in \Pic0$ we have $h^0
( E \otimes \xi )= h^1 (E \otimes \xi)$. By base change, the
condition that $E$ is not a base point is equivalent to the fact
that $R^0 \Phi_P E$ is supported on a proper subset of $\Pic0$, but
since $R^0 \Phi_P E$ is torsion-free\footnote{For example embed $E$
in $E(D)$ for some divisor $D$ on $X$ of very large degree and apply
$\Phi_P$ to the inclusion.}, this is equivalent to requiring it to
be $0$.
\end{proof}

\begin{remark}
Via an Abel-Jacobi embedding of $X$ in its Jacobian $J(X)$, the
functor $\R \Phi_P$ is the same as the Fourier-Mukai transform $\R
\widehat \SS$ on $J(X)$ applied to objects supported on $X$. Hence
in the statement here we might as well talk about $\R \widehat \SS$
instead of $\R \Phi_P$.
\end{remark}

Fix now any polarization $\Theta$ on $J(X)$, and consider for any
$m$ the Fourier-Mukai transform
$$E_m : = \R \Psi_P \OO_{J(X)} (-m\Theta) [g] = \widehat{\OO_{J(X)} (-m\Theta)},$$
which by base change is a vector bundle on $X$, usually called a
\emph{Raynaud bundle}. Using Lemma \ref{base1}, together with the
general statements of Corollary \ref{k=0} and Lemma
\ref{cohomology}, we obtain the following criterion for detecting
base points for the determinant line bundle. The equivalence of (1)
and (3) has already been noted by Hein \cite{hein}, Theorem 2.5.
With a more careful study, he gives an effective bound for $m$ in
(3) (cf. \emph{loc.cit}, Theorem 3.7).

\begin{corollary}\label{base_criterion_1}
Let $X$ be a smooth projective curve of genus $g \ge 2$ and $E$ a
vector bundle in $SU_X(r, L)$, with $L \in {\rm Pic}^{g-1}(X)$. Then
the following are equivalent:\\
(1) $E$ is not a base point for the linear series $|\LL|$.\\
(2) $H^i (E \otimes E_m) = 0$, for all $i >0$ and all $m \gg0$.\\
(3) $H^0 (E \otimes E_m^\vee) = 0$ for all $m \gg0$.
\end{corollary}

The Strange Duality conjecture, proved recently by Belkale
\cite{belkale} and Marian-Oprea \cite{mo}, allows for an extension
of this in the most general setting. Let $SU_X(r,L)$ be as above,
with $L \in {\rm Pic}^d(X)$. Let $h := (r,d)$, and $r_0 := r/h$ and
$d_0 := d/h$. A general bundle $F$ in the moduli space $U_X(kr_0,
k(r_0 (g-1) -d_0))$ (with arbitrary determinant) gives a generalized
theta divisor
$$\Theta_F := \{ E ~|~h^0 (E \otimes F) \neq 0\} \subset SU_X(r,L)$$
which belongs to the linear series $|\LL^k|$ (cf. \cite{dn}). The
Strange Duality is equivalent to the fact that the divisors
$\Theta_F$ span this linear series as $F$ varies. It is well known
that there exists as cover $M$ of $U_X(kr_0, k(r_0 (g-1)
-d_0))$, \'etale over the stable locus, such that there is a universal bundle 
$\E$ on $X\times M$.
We consider the Fourier-Mukai correspondence $\R\Phi_{\E}: \D(X)
\rightarrow \D(M)$.\footnote{For the statements we are interested
in, this is a good as thinking that $M$ is $U_X(kr_0, k(r_0 (g-1)
-d_0))$ itself, with the technical problem that as soon as $k \ge 2$
this moduli space will definitely not be fine.} We thus obtain as
before:

\begin{lemma}\label{base2}
If $E$ corresponds to a point in $SU_X(r,L)$, the following are
equivalent:\\
(1) $E$ is not a base point for $|\LL^k|$.\\
(2)$E$ has a theta divisor $\Theta_E := \{F ~|~h^0 (E \otimes F)
\neq 0\} \subset U_X(kr_0, k(r_0 (g-1) -d_0))$.\\
(3) $E$ is a $GV$-sheaf with respect to any universal bundle $\E$ on
$X \times M$.\\
(4) $R^0 \Phi_{\E} E = 0$, i.e. $E$ satisfies $WIT_1$ with respect
to $\E$.
\end{lemma}

We fix any generalized theta divisor on $U_X(kr_0, k(r_0 (g-1)
-d_0))$, and denote abusively by $\Theta$ its pullback to $M$. We
consider for $m \gg0$ the Fourier-Mukai transform
$$E_m^k : = \R \Psi_{\E} \OO_M (-m\Theta) [\dim M] = \widehat {\OO_M(-m\Theta)},$$
which is a vector bundle on $X$ generalizing Raynaud's bundles
coming from the Jacobian. (This bundle can be constructed also as the push-forward of a Raynaud bundle
on a spectral cover of $C$ associated to the moduli space $U_X(kr_0, k(r_0 (g-1)
-d_0))$ as in \cite{bnr}.) As above, this provides the promised
extension of Corollary \ref{base_criterion_1}.

\begin{corollary}\label{base_criterion_2}
Let $X$ be a smooth projective curve of genus $g \ge 2$ and $E$ a
vector bundle in $SU_X(r, L)$, with $L \in {\rm Pic}^{d}(X)$. Then
the following are equivalent:\\
(1) $E$ is not a base point for the linear series $|\LL^k|$.\\
(2) $H^i (E \otimes E_m^k) = 0$, for all $i >0$ and all $m \gg0$.\\
(3) $H^0 (E \otimes {E_m^k}^\vee) = 0$ for all $m \gg0$.
\end{corollary}

\noindent {\bf Relative moduli of sheaves on threefold Calabi-Yau
fibrations.} In theory one can study generic vanishing statements
for any setting of the type: $X$ is smooth projective, $M$ is a fine
moduli space of objects over $X$, and $\E$ is a universal object on
$X\times M$ inducing the functor $\Phi _{\E}$. In practice, the main
difficulty to be overcome is a good understanding of the vector
bundles  $\widehat{A^{-1}} = \R \Psi_{\E} (A^{-1})[{\rm dim}~M]$,
with $A$ a very positive line bundle on $M$, on a case by case
basis. Very few concrete examples seem to be known beyond the case
of abelian varieties.\footnote{Besides curves or surfaces, where the
vanishing of cohomology groups of appropriate semistable sheaves can
usually be tested by hand -- cf. for example \S4, listing various
examples of Yoshioka.} We would like to raise as a general problem
to describe the structure of these vector bundles, given a specific
moduli space.\footnote{For example, on a curve $X$, in the notation
of the previous section the bundles $E_m$ are well understood up to
isogeny, due to Mukai's results on abelian varieties. How about the
bundles $E_m^k$?}

Here we give only a rather naive example of such a result for a
Fourier-Mukai functor associated to threefolds with abelian or $K3$
fibration, considered first by Bridgeland and Maciocia in \cite{bm}
(cf. also \cite{bridgeland}). This works under special numerical
hypotheses, based on results of Yoshioka. The interested reader can
prove similarly an analogous result in the case of elliptic
threefolds.

Recall that a \emph{Calabi-Yau fibration} is a morphism $\pi:
X\rightarrow S$ of smooth projective varieties, with connected
fibers, such that $K_X\cdot C=0$ for all curves $C$ contained in
fibers of $\pi$. If it is of relative dimension at most two, then it
is an elliptic, abelian surface, or $K3$-fibration (in the sense
that the nonsingular fibers are of this type). Say $\pi$ is flat,
and consider a polarization $H$ on $X$, and $Y$ an irreducible
component of the relative moduli space $M^{H, P}(X/S)$ of sheaves on
$X$ (over $S$), semistable with respect to $H$, and with fixed
Hilbert polynomial $P$. The choice of $P$ induces on every smooth
fiber $X_s$ invariants which are equivalent to the choice of a Mukai
vector $v \in H^{ev}(X_s, \ZZ)$ as in \S4(5). Assuming that $Y$ is
also a threefold, and fine, Bridgeland and Maciocia (cf. \cite{bm},
Theorem 1.2) proved that it is smooth, and the induced morphism
$\hat \pi: Y \rightarrow S$ is a Calabi-Yau fibration of the same
type as $\pi$. In addition, if $\E$ is a universal sheaf on $X\times
Y$, then the Fourier-Mukai functor $\R \Phi_{\E}: \D(X) \rightarrow
\D(Y)$ is an equivalence of derived categories. We consider the
following condition:

\noindent ($\star$) ~ For each $s\in S$ such that $X_s$ is smooth,
the Mukai vector $v$ is primitive and isotropic, and the structure
sheaf $\OO_{X_s}$ satisfies $WIT_2$ with respect to the induced $\R
\Phi_{\E_s}: \D(X_s ) \rightarrow \D(Y_s)$.

\begin{example}
Papers of Yoshioka (e.g. \cite{yoshioka2}, \cite{yoshioka3}) contain
plenty of examples of surfaces where condition ($\star$) is
satisfied. For some precise ones, both abelian and $K3$, cf. \S4(5).
Note that in all the cases we know, we have ${\rm Pic}(X_s) \cong
\ZZ$.
\end{example}

If the first half of ($\star$) is satisfied, it is proved in
\cite{bm} \S7 that the moduli space $M^{H, P}(X/S)$ does have a fine
component $Y$ which is a threefold, so the above applies. For
simplicity we assume in the next statement that all the fibers of
$\pi$ are smooth, but please note Remark \ref{singular_fibers},
which explains that the result can be made more general.

\begin{proposition/example}\label{cy}
Let $X$ be a smooth projective threefold with a smooth Calabi-Yau
fibration $\pi: X \rightarrow S$ of relative dimension two. Let $H$
be a polarization on $X$ and $P$ a Hilbert polynomial, and assume
that condition ($\star$) is satisfied. Consider a fine
three-dimensional moduli space component $Y \subset M^{H,P} (X/S)$,
and let $\E$ be a universal sheaf on $X\times Y$. Then $\omega_X $
is a $GV_{-1}$-sheaf with respect to $\E$. In particular
$$H^i (X, \omega_X\otimes E) = 0, ~\forall ~i > 1, ~\forall~ E\in Y {\rm~general}.$$
\end{proposition/example}
\begin{proof}
In order to prove that $\omega_X$ is $GV_{-1}$ with respect to $\E$, it
is enough to check condition (2) in Theorem \ref{general_f}. Given a
very positive line bundle $A$ on $Y$ we want
\begin{equation}\label{condition}
H^i (\omega_X \otimes \widehat {A^{-1}} ) = 0, ~{\rm for~all~} i >
1,
\end{equation}
where the Fourier transform is with respect to $\Psi_{\E}$.

We use the facts established in \cite{bm}: the moduli space
$M^{H,P}(X/S)$ restricts for each $s\in S$ to the corresponding
moduli spaces of sheaves with Mukai vector $v$ on $X_s$, stable with
respect to the polarization $H_{|X_s}$. The functor $\R\Phi_{\E}$ is
a relative Fourier-Mukai functor, which induces the respective
fiberwise functors $\R\Phi_{\E_s}: \D(X_s) \rightarrow \D(Y_s)$.
With our choice of polarization, each $Y_s$ is a fine moduli space
of sheaves on $X_s$, of the same dimension and in fact the same type
as $X_s$ (cf. \cite{bm}, \S7.1, using Mukai's results).\footnote{In
fact each $\R \Phi_{\E_s}$ is an equivalence of derived categories.}

We check condition (\ref{condition}) by using the Leray spectral
sequence for $\pi: X \rightarrow S$, namely
$$E^{i,j}_2 := H^i (S, R^j \pi_* (\omega_X \otimes \widehat {A^{-1}})) \Rightarrow
H^{i+j} (X, \omega_X \otimes \widehat {A^{-1}}).$$ For every $s \in
S$ we have $(\omega_X \otimes \widehat {A^{-1}})_{|X_s} \cong
\widehat {A_s^{-1}}$, where we denote $A_s : = A_{|Y_s}$, and the
transform on the right hand side is taken with respect to $\R
\Psi_{\E_s}$.

By assumption we have that $\OO_{X_s}$ satisfies $WIT_2$ with
respect to $\R\Phi_{\E_s}$. We can then use Corollary \ref{k=0} in a
different direction ((3) $\Rightarrow$ (2)), to deduce that $H^i
(X_s, \widehat {A_s^{-1}}) = 0$. Note that ${\omega_X} _{|X_s} \cong
\omega_{X_s} \cong \OO_{X_s}$. Since $\pi$ is smooth, we obtain by
base-change that $R^j \pi_* (\omega_X \otimes \widehat {A^{-1}}) =
0$ for all $j \ge 1$. This immediately gives that $E^{i,j}_2 = 0$
for $i\ge 2$ and all $j$, for $i = 1$ and $j \ge 1$, and also for $i
= 0$ and $j =2$.  Thus the spectral sequence provides
$$H^{i} (X, \omega_X \otimes \widehat {A^{-1}}) = 0 , ~{\rm for~} i = 2,3.$$
This is precisely (\ref{condition}), and we get that $\omega_X$ is
$GV_{-1}$, or equivalently by Theorem \ref{f}, that  $R^i \Phi \OO_X =
0$ for $i < 2$.
\end{proof}

\begin{remark}\label{singular_fibers}
The same proof works in fact if we don't necessarily assume that
$\pi$ is smooth, but only that it is flat (a necessary assumption),
plus the slightly technical condition $R^2 \pi_* (\omega_X \otimes
\widehat {A^{-1}}) = 0$ for $A$ sufficiently positive on
$Y$.\footnote{It should however be true that this holds most of the
time, though admittedly we have not checked -- at least in the case
of abelian fibrations, after a relative isogeny one expects to
reduce the calculation to a sheaf of the form $R^2 \pi_* (\omega_X
\otimes M)$ with $M$ semiample, which would be $0$ by Koll\'ar's
torsion-freeness theorem.} Indeed, since there is only a finite
number of singular fibers, we obtain by base-change that $R^j \pi_*
(\omega_X \otimes \widehat {A^{-1}})$ is supported on at most a
finite set, for $j = 1,2$, and it is of course $0$ for $j \ge 3$.
This immediately gives that $E^{i,j}_2 = 0$ for $i\ge 2$ and all
$j$, and also for $i = 1$ and $j \ge 1$. The only term which may
cause trouble is $E^{0,2}_2 = H^0 (S, R^2 \pi_* (\omega_X \otimes
\widehat {A^{-1}}))$, and for its vanishing we have to use the
assumption above.
\end{remark}

\providecommand{\bysame}{\leavevmode\hbox
to3em{\hrulefill}\thinspace}

\end{document}